\documentclass[12pt,a4paper]{amsart}
\usepackage{a4wide}
\usepackage[T1]{fontenc}
\usepackage[utf8]{inputenc}
\usepackage{lmodern}
\usepackage{microtype}
\usepackage{amsmath,amssymb,amsthm,mathtools}
\usepackage{enumitem}
\usepackage{url}
\usepackage{advdate}
\usepackage[colorlinks=true,linkcolor=blue,citecolor=blue,urlcolor=blue]{hyperref}
\hypersetup{pdftitle={The e-Property for Quantum Markov Semigroups},pdfauthor={Tomasz Kania and Mariusz Niwinski}}
\numberwithin{equation}{section}

\title[The $e$-property for quantum Markov semigroups]{The $e$-Property for Quantum Markov Semigroups}
\author[T.~Kania]{Tomasz Kania}
\address[T.~Kania]{Mathematical Institute\\Czech Academy of Sciences\\\v Zitn\'a 25 \\115 67 Praha 1\\Czech Republic  and  Institute of Mathematics and Computer Science\\ Jagiellonian University\\ {\L}ojasiewicza 6, 30-348 Krak\'{o}w, Poland}
\email{kania@math.cas.cz, tomasz.marcin.kania@gmail.com}
\thanks{RVO: 67985840.}
\author[M.~Niwi\'nski]{Mariusz Niwi\'nski}
\address[M.~Niwi\'nski]{Doctoral School of Exact and Natural Sciences\\Faculty of Mathematics and Computer Science\\Institute of Mathematics, Jagiellonian University\\{\L}ojasiewicza 6, 30-348 Krak\'ow, Poland}
\email{mariusz.niwinski@doctoral.uj.edu.pl}
\thanks{The second-named author acknowledges with thanks support from NCN SONATA BIS~13 (2023/50/E/ST1/00067).}
\date{October 24, 2026}

\newtheorem{theorem}{Theorem}[section]
\newtheorem{proposition}[theorem]{Proposition}
\newtheorem{corollary}[theorem]{Corollary}
\newtheorem{lemma}[theorem]{Lemma}
\theoremstyle{definition}
\newtheorem{definition}[theorem]{Definition}
\newtheorem{remark}[theorem]{Remark}
\newtheorem*{remark*}{Remark}
\newtheorem{example}[theorem]{Example}
\newtheorem{problem}[theorem]{Problem}
\newtheorem*{theoremA}{Theorem A}
\newtheorem*{theoremB}{Theorem B}
\newtheorem*{theoremC}{Theorem C}
\newtheorem*{theoremD}{Theorem D}

\newcommand{\A}{\mathcal A}
\newcommand{\B}{\mathcal B}
\newcommand{\C}{\mathbb C}
\newcommand{\D}{\mathcal D}
\newcommand{\M}{\mathcal M}
\newcommand{\Mstar}{\mathcal M_*}
\newcommand{\N}{\mathbb N}
\newcommand{\Nc}{\mathcal N}
\newcommand{\R}{\mathbb R}
\newcommand{\Rr}{\mathcal R}
\newcommand{\Sst}{\mathcal S}
\newcommand{\Ldom}{\mathcal L}

\newcommand{\T}{\mathcal T}
\newcommand{\id}{\mathrm{id}}
\newcommand{\AP}{\operatorname{AP}}
\newcommand{\Inv}{\operatorname{Inv}}
\newcommand{\MD}{\operatorname{MD}}
\newcommand{\pn}{\mathrm{pn}}
\newcommand{\eDef}{\mathfrak e}
\newcommand{\chiT}{\boldsymbol\chi}
\DeclareMathOperator{\dom}{dom}
\DeclareMathOperator{\Lip}{Lip}
\DeclareMathOperator{\Fix}{Fix}
\DeclareMathOperator{\diam}{diam}
\DeclareMathOperator{\osc}{osc}

\begin{document}
\begin{abstract}
Let $(\T_t)_{t\geqslant0}$ be a quantum Markov semigroup on a von Neumann algebra and let $\A$ be an invariant Feller $C^*$-algebra.  We formulate the $e$-property on the full state space and show that its quantitative defect is equivalent, within an optimal factor two, to the Hausdorff measure of non-compactness of observable orbits.  We develop Ces\`aro and Lip-norm criteria and a generalised Doeblin minorisation by a state-space contracting positive map, including finite-dimensional and measure-prepare factorisations.  When every observable has the $e$-property, the Jacobs--de Leeuw--Glicksberg decomposition identifies the reversible Choi--Effros algebra and the norm-stable sector.  Expected and decoherence-free subalgebras transfer the property under norm decay of the complementary dynamics, including infinite-dimensional ultracontractive and finite-index Pimsner--Popa regimes. Finally, we establish the full comparison hierarchy among ordinary, eventual, Ces\`aro and pointwise eventual equicontinuity notions, reinforced by matrix Jamison and homoclinic flow counterexamples. Projection-face criteria are treated in a companion paper.
\end{abstract}
\maketitle

\section{Introduction}

For a classical Markov--Feller semigroup $(P_t)_{t\geqslant0}$ on a Polish metric space $X$, the $e$-property asks that $\{P_tf:t\geqslant0\}$ be equicontinuous for every bounded Lipschitz function $f$.  It is a uniform-in-time continuity principle for expectations and enters invariant-measure theory, asymptotic stability, lower-bound arguments and iterated function systems.  The results of Czapla on Markov--Feller $e$-chains and partially equicontinuous operators, the work of Czapla--Horbacz on iterated systems, and later comparisons of ordinary, eventual and averaged equicontinuity are among the direct motivations for this paper \cite{Czapla2012,Czapla2018,CzaplaHorbacz2014,CzaplaHilleHorbaczWojewodka2020,HilleSzarekZiemlanska2017,HilleSzarekWormZiemlanska2021,KukulskiWojewodka2024,LiuLiu2024,Szarek2000,SzarekWorm2012}.  The classical $e$-property viewpoint is therefore the organising principle of the paper.

The quantum framework has two layers.  A quantum Markov semigroup acts normally and weak-$^*$ continuously on a von Neumann algebra $\M$.  We simultaneously retain an invariant unital $C^*$-subalgebra $\A\subseteq\M$ on which every observable orbit is norm continuous.  The ambient von Neumann algebra carries the standard quantum dynamics and its predual, whereas the Feller algebra $\A$ carries compact initial orbit segments, the weak-$^*$ compact state space and the norm geometry of the $e$-property.  Thus the $C^*$-algebra is not a replacement for the von Neumann setting; it is the norm-continuous layer already natural in quantum Feller theory \cite{CiprianiFagnolaLindsay2000,CarboneFagnola2003}.

For $a=a^*\in\A$, we study the expectation functions
\[
\Sst(\A)\ni\varphi\longmapsto\varphi(\T_t(a)) \in \R,
\qquad t\geqslant0.
\]
Their weak-$^*$ equicontinuity is the full-state-space $e$-property.  In a compact commutative system it agrees with the ordinary pointwise property, while in a non-commutative algebra the whole state space provides an intrinsic compact phase space.

\begin{theoremA}
For every $a=a^*\in\A$, the equicontinuity defect $\eDef_{\T}(a)$ and the Hausdorff measure of non-compactness $\chiT_{\T}(a)$ satisfy
\[
\frac12\eDef_{\T}(a)
\leqslant\chiT_{\T}(a)
\leqslant\eDef_{\T}(a).
\]
Both quantities are unchanged by deleting a bounded initial time interval and are Lipschitz under perturbations of $a$.  Consequently the $e$-property, its eventual version and norm relative compactness of the orbit are equivalent.
\end{theoremA}

The factor two is optimal for bounded families on compact spaces.  We also study Ces\`aro averages: the ordinary $e$-property implies the averaged one, and norm convergence of ergodic averages yields the Ces\`aro $e$-property.  Compact quantum metrics provide another route through uniform or strict Lip-norm estimates.

The second theme is global mixing.  The classical Doeblin condition dominates one time step by a reference probability. Here the minorising map itself may already contract the state space.

\begin{theoremB}
Let $T\colon\A\to\A$ be unital and positive.  Suppose that
\[
T(b)\geqslant\lambda\Psi(b)\qquad(b\in\A_+),
\]
where $\Psi$ is unital and positive, $\lambda\in(0,1]$, and
\[
\|\varphi\circ\Psi-\psi\circ\Psi\|
\leqslant\beta\|\varphi-\psi\|
\qquad(\varphi,\psi\in\Sst(\A))
\]
for some $\beta<1$.  Then the dual of $T$ contracts the state space with coefficient at most $1-\lambda(1-\beta)$.  For a time slice of the semigroup this gives a unique invariant state, exponential dual-norm convergence and the $e$-property on all self-adjoint observables.
\end{theoremB}

Positive return profiles, composable minorizations, positivity-improving finite-dimensional channels and measure-prepare maps are used to produce rank-one lower bounds. This global theory is distinct from the local corner mechanism developed in the companion paper \cite{KaniaNiwinskiProjectionFaces}.

When every observable has the $e$-property, the semigroup is strongly almost periodic.  The Jacobs--de Leeuw--Glicksberg architecture is established Banach-space and quantum semigroup theory \cite{Hellmich2011,BatkaiGrohKunszentiSchreiber2012,RundeViselter2014}; we use it as the structural meaning of the $e$-property rather than presenting the abstract decomposition as new.

\begin{theoremC}
If every orbit in $\A$ is relatively compact in norm, there is a unital completely positive projection $P$ commuting with the semigroup such that
\[
\A=P(\A)\oplus\ker P,
\qquad
\|\T_t(a)\|\longrightarrow0\quad(a\in\ker P),
\]
and the point-norm closure of the restricted dynamics on $P(\A)$ is a compact abelian group of complete order automorphisms.  The range is a $C^*$-algebra for the Choi--Effros product.  Its largest ambient $C^*$-subalgebra is $P(\A)\cap\MD(P)$.
\end{theoremC}

We give norm spectral synthesis, describe invariant states and provide infinite-dimensional examples where the Choi--Effros and ambient products differ.

\begin{theoremD}
Let $E\colon\A\to\B$ be a conditional expectation onto an invariant unital $C^*$-subalgebra and put
\[
\eta(a)=\limsup_{t\to\infty}\|\T_t(a-E(a))\|.
\]
Then the orbit non-compactness measures of $a$ and $E(a)$ differ by at most $\eta(a)$, while their $e$-defects differ by at most $2\eta(a)$.  If $\eta(a)=0$, they have the same norm omega-limit set and the same norm convergence behaviour.
\end{theoremD}

This yields lifting results for Cartan masas and finite-dimensional spectral gaps.  It also treats decoherence-free subalgebras.  Throughout the paper, the $e$-property remains the organising thread.  We first treat its state-space formulation and quantitative defect, then its Ces\`aro and reversible consequences, its metric and Doeblin criteria, and finally transfer from expected subalgebras.

\section{The von Neumann--Feller framework}\label{sec:framework}

\begin{definition}
Let $\M$ be a von Neumann algebra. A \emph{quantum Markov semigroup} on $\M$ is a family $(\T_t)_{t\geqslant 0}$ of normal, unital, completely positive maps on $\M$ such that
\begin{enumerate}[label=\upshape(\roman*)]
    \item $\T_0 = \id$,
    \item $\T_{s+t} = \T_s \T_t$ for all $s,t\geqslant 0$,
    \item for every $x\in\M$ and every $\omega\in\Mstar$, the map $t\mapsto \omega(\T_t(x))$ is continuous on $[0,\infty)$.
\end{enumerate}
\end{definition}

\begin{remark}
The predual maps $\T_{t,*}\colon \Mstar\to\Mstar$, defined by
\[
(\T_{t,*}\omega)(x) = \omega(\T_t(x)), \qquad \omega\in\Mstar,\ x\in\M,
\]
preserve the normal state space of $\M$.
\end{remark}

\begin{definition}
Let $(\T_t)_{t\geqslant 0}$ be a quantum Markov semigroup on $\M$, and let $\A\subseteq\M$ be a unital $C^*$-subalgebra. We say that $\A$ is a \emph{Feller algebra} for $(\T_t)$ if
\begin{enumerate}[label=\upshape(\roman*)]
    \item $\T_t(\A)\subseteq\A$ for every $t\geqslant 0$,
    \item for every $a\in\A$, the map $t\mapsto \T_t(a)$ is norm continuous from $[0,\infty)$ into $\A$.
\end{enumerate}
In this case we call $(\T_t)_{t\geqslant 0}$ a \emph{quantum Markov--Feller semigroup} on $\A$.
\end{definition}

\begin{remark}\label{rem:hybrid-setting}
The von Neumann algebra and the Feller algebra have complementary roles.  Normality and weak-$^*$ continuity belong to the ambient quantum Markov dynamics on $\M$, whereas norm-continuous orbit segments, the compact state space and the norm compactness underlying the $e$-property live on $\A$.  Thus the present framework does not replace von Neumann algebras by abstract $C^*$-algebras.  Some order estimates below use only the strongly continuous unital positive semigroup on $\A$; this extra generality is stated explicitly.
\end{remark}

\begin{remark}
We keep $\A$ unital in order to have a compact state space and to parallel the bounded-observable picture. In the locally compact commutative case one may start from $C_0(X)$ and pass to its multiplier algebra $M(C_0(X))=C_b(X)$. More generally, non-degenerate positive maps on a $C_0$-type algebra extend, under the usual strict-continuity hypotheses, to multiplier algebras, and the state-space definitions below may then be applied on $M(\A)$.
\end{remark}

Let $\Sst(\A)$ denote the state space of $\A$, endowed with its weak-$^*$ topology. We write $\T_t^*\colon\A^*\to\A^*$ for the Banach-space dual of $\T_t|_{\A}$, so that $\T_t^*(\varphi)=\varphi\circ\T_t$; the notation $\T_{t,*}$ is reserved for the predual action on $\Mstar$. Throughout the paper, $\Ldom\subseteq\A_{\mathrm{sa}}$ denotes a unital real linear test space. Whenever norm density is needed, it is stated explicitly. We use $\Ldom$ as the class of test observables.

\begin{definition}
Let $(\T_t)_{t\geqslant 0}$ be a quantum Markov semigroup on $\M$ and let $\A\subseteq\M$ be a Feller algebra.
\begin{enumerate}[label=\upshape(\alph*)]
    \item We say that $(\T_t)$ has the \emph{$e$-property at $\varphi\in\Sst(\A)$ relative to $\Ldom$} if for every $a\in\Ldom$ and every $\varepsilon>0$ there exists a weak-$^*$ neighbourhood $U$ of $\varphi$ in $\Sst(\A)$ such that
    \[
    \sup_{\psi\in U}\sup_{t\geqslant 0}\bigl|\psi(\T_t(a)) - \varphi(\T_t(a))\bigr| < \varepsilon.
    \]
    \item We say that $(\T_t)$ has the \emph{eventual $e$-property at $\varphi\in\Sst(\A)$ relative to $\Ldom$} if for every $a\in\Ldom$ and every $\varepsilon>0$ there exist a weak-$^*$ neighbourhood $U$ of $\varphi$ in $\Sst(\A)$ and a time $t_0\geqslant 0$ such that
    \[
    \sup_{\psi\in U}\sup_{t\geqslant t_0}\bigl|\psi(\T_t(a)) - \varphi(\T_t(a))\bigr| < \varepsilon.
    \]
\end{enumerate}
If the corresponding property holds at every $\varphi\in\Sst(\A)$, we simply say that $(\T_t)$ has the $e$-property, respectively the eventual $e$-property, on $\Ldom$.
\end{definition}

\begin{remark}\label{rem:norm-vs-weak}
The dual maps $\T_t^*$ are contractive in the norm of $\A^*$, as is the dual of every unital positive map. Therefore, for every $t\geqslant 0$ and every $\varphi,\psi\in\Sst(\A)$, one has
\[
\|\T_t^*(\varphi)-\T_t^*(\psi)\|
=
\|\varphi\circ\T_t-\psi\circ\T_t\|
\leqslant
\|\varphi-\psi\|.
\]
Consequently the family of maps
\[
\Sst(\A)\ni\varphi\longmapsto \varphi\circ\T_t \in \Sst(\A)
\]
is always uniformly equicontinuous with respect to the \emph{norm} topology on the state space.

This is not the $e$-property. In infinite dimensions the norm topology on $\Sst(\A)$ is far too strong and, in particular, is not compact. The classical $e$-property, and our non-commutative formulation, concern instead the much coarser weak-$^*$ topology, or equivalently a metric topology such as the Monge--Kantorovich topology when a Lip-norm is available. This mirrors the familiar probabilistic distinction between total variation --- always contractive for Markov operators --- and weak convergence of measures, where equicontinuity in Wasserstein or bounded-Lipschitz type distances becomes a genuinely rigid additional property.
\end{remark}

\begin{remark}\label{rem:bures-trace}
The preceding remark also explains why trace-norm or fidelity-based metrics do not by themselves capture the phenomenon studied here. On normal states of a von Neumann algebra, the predual maps of a quantum Markov semigroup are contractions in trace norm. In finite dimensions, and more generally in the usual normal-state settings where fidelity is available, the data-processing inequality gives
\[
F(\varphi\circ \T_t,\psi\circ \T_t)\geqslant F(\varphi,\psi),
\]
and hence the corresponding Bures distance is non-increasing under the dual dynamics. Thus trace distance and Bures distance give automatic global equicontinuity for every quantum Markov semigroup. The non-trivial issue is instead equicontinuity for the weak-$^*$ topology, or for Monge--Kantorovich metrics that induce it. This is why the present paper works with state-space weak-$^*$ equicontinuity rather than with the standard information-theoretic contraction metrics.
\end{remark}

\begin{remark}
The classical literature distinguishes several eventual or averaged equicontinuity notions; see, for example, \cite{HilleSzarekWormZiemlanska2021,Czapla2018,KukulskiWojewodka2024,LiuLiu2024}. Besides the ordinary $e$-property and the Ces\`aro $e$-property, one encounters at least the following variants:
\begin{enumerate}[label=\upshape(\roman*),leftmargin=2.6em]
    \item the \emph{fixed-threshold eventual $e$-property}, where one asks for a single time $\tau$ such that
    \[
    \lim_{x\to z}\sup_{t\geqslant \tau}|P_t f(x)-P_t f(z)|=0;
    \]
    \item the weaker \emph{$\varepsilon$-dependent eventual $e$-property}, where the threshold time may depend on the tolerance $\varepsilon$;
    \item \emph{eventual continuity}, where one only asks that the large-time tails become equicontinuous in the limit, without requiring one fixed threshold time to work uniformly;
    \item the \emph{Ces\`aro $e$-property}, obtained by replacing the orbit by its time averages.
\end{enumerate}
Our definition is the version in \upshape(ii): for each observable and each tolerance $\varepsilon$, one weak-$^*$ neighbourhood works uniformly for all $t\geqslant t_0(\varepsilon)$.

The relations most relevant for the present paper may be summarised schematically as
\[
\text{fixed-threshold eventual }e
\Longrightarrow
\text{eventual }e\text{ used here}
\Longrightarrow
\text{eventual continuity},
\]
\[
\text{eventual }e\text{ used here}
\xrightarrow{\ \text{strong stochastic continuity}\ }
e
\Longrightarrow
\text{Ces\`aro }e.
\]
The second implication is classical for strongly stochastically continuous Markov--Feller semigroups \cite{KukulskiWojewodka2024}; in the present non-commutative framework it is replaced by Theorem~\ref{thm:eventual-implies-e}, where norm continuity on the Feller algebra plays the corresponding rôle. For Markov operators, several equicontinuity notions become equivalent under additional Schur-like hypotheses on spaces of measures \cite{HilleSzarekWormZiemlanska2021}, while the relation with eventual continuity is analysed further in \cite{LiuLiu2024}. In general the reverse arrows fail, as shown by examples in \cite{HilleSzarekZiemlanska2017,Czapla2018,KukulskiWojewodka2024}. Historically, the terminology is also connected with the theory of $e$-chains in the sense of Meyn and Tweedie \cite{MeynTweedie2009}, although the present paper follows the more recent $e$-property language used in the Markov--Feller literature. We use the version above because it is the cleanest one to transfer to the non-commutative state-space setting and it already suffices for the compactness argument below.
\end{remark}

\begin{remark}
This definition is deliberately formulated on the full state space of $\A$, not on pure states and not on normal states of $\M$. In the commutative case points are evaluation states; in the genuinely non-commutative case the whole state space is the robust substitute for the classical phase space.
\end{remark}

\begin{remark}
Let $X$ be a Polish metric space, let $\A=C_b(X,\C)$, let $\Ldom=\Lip_b(X,\R)$ viewed simply as a unital real test space, and let $(P_t)_{t\geqslant 0}$ be a classical Markov--Feller semigroup. Then for every $x\in X$ the evaluation functional $\delta_x\in\Sst(\A)$ satisfies
\[
\delta_y(P_t f)-\delta_x(P_t f)=P_t f(y)-P_t f(x).
\]
Hence the classical $e$-property is exactly the restriction of the present state-space definition to the embedded copy $\{\delta_x : x\in X\}\subseteq\Sst(\A)$.

Two remarks are worth making. First, for a general Polish metric space $X$ the space $\Lip_b(X,\R)$ need not be norm dense in $C_b(X,\R)$, so in the present paper it should be viewed as a test space rather than as a dense subspace. If one wants norm density, one may either work on a compact metric space $X$, where $\Lip(X,\R)$ is uniformly dense in $C(X,\R)$, or, when available, work on an invariant commutative Feller subalgebra containing the bounded Lipschitz functions. Secondly, the full state space of $C_b(X)$ is larger than the copy of $X$ itself. The literal classical phase space is therefore recovered, exactly as intended, by restricting to the Dirac states.
\end{remark}

\section{From the eventual \texorpdfstring{$e$}{e}-property to the \texorpdfstring{$e$}{e}-property}

The key compactness input is completely elementary.

\begin{lemma}\label{lem:compact-uniform}
Let $\A$ be a $C^*$-algebra, let $K\subseteq \A$ be norm compact, and let $(\varphi_i)_i$ be a net in $\Sst(\A)$ converging to $\varphi$ in the weak-$^*$ topology. Then
\[
\sup_{a\in K}\bigl|\varphi_i(a)-\varphi(a)\bigr|\longrightarrow 0.
\]
\end{lemma}

\begin{proof}
Fix $\varepsilon>0$. Choose $a_1,\dots,a_n\in K$ such that every $a\in K$ satisfies $\min_k\|a-a_k\|<\varepsilon/3$. Since $\varphi_i\to\varphi$ weak-$^*$, there exists $i_0$ such that
\[
\bigl|\varphi_i(a_k)-\varphi(a_k)\bigr|<\varepsilon/3,\qquad i\geqslant i_0,\ 1\leqslant k\leqslant n.
\]
Take $i\geqslant i_0$ and $a\in K$, and choose $k$ with $\|a-a_k\|<\varepsilon/3$. Because states are contractive,
\begin{align*}
\bigl|\varphi_i(a)-\varphi(a)\bigr|
&\leqslant \bigl|\varphi_i(a-a_k)\bigr| + \bigl|\varphi_i(a_k)-\varphi(a_k)\bigr| + \bigl|\varphi(a_k-a)\bigr| \\
&\leqslant 2\|a-a_k\| + \varepsilon/3 < \varepsilon.
\end{align*}
This proves the claim.
\end{proof}

\begin{theorem}\label{thm:eventual-implies-e}
Let $(\T_t)_{t\geqslant 0}$ be a quantum Markov semigroup on a von Neumann algebra $\M$, and let $\A\subseteq\M$ be a Feller algebra. If $(\T_t)$ has the eventual $e$-property on $\Ldom$, then it has the $e$-property on $\Ldom$.
\end{theorem}

\begin{proof}
Fix $\varphi\in\Sst(\A)$, $a\in\Ldom$, and $\varepsilon>0$. By the eventual $e$-property there exist a weak-$^*$ neighbourhood $U_1$ of $\varphi$ and a time $t_0\geqslant 0$ such that
\[
\sup_{\psi\in U_1}\sup_{t\geqslant t_0}\bigl|\psi(\T_t(a)) - \varphi(\T_t(a))\bigr| < \varepsilon/2.
\]
Since $\A$ is a Feller algebra, the map $[0,t_0]\ni t\mapsto \T_t(a)\in\A$ is norm continuous. Consequently,
\[
K:=\{\T_t(a): 0\leqslant t\leqslant t_0\}
\]
is norm compact. We claim that there exists a weak-$^*$ neighbourhood $U_2$ of $\varphi$ such that
\[
\sup_{\psi\in U_2}\sup_{b\in K}\bigl|\psi(b)-\varphi(b)\bigr| < \varepsilon/2.
\]
Indeed, otherwise one could find a net $(\psi_i)$ in $\Sst(\A)$ with $\psi_i\to\varphi$ weak-$^*$ and
\[
\sup_{b\in K}\bigl|\psi_i(b)-\varphi(b)\bigr|\geqslant \varepsilon/2
\]
for all $i$, contradicting Lemma~\ref{lem:compact-uniform}. For $U:=U_1\cap U_2$ we obtain, for every $\psi\in U$,
\begin{align*}
\sup_{t\geqslant 0}\bigl|\psi(\T_t(a)) - \varphi(\T_t(a))\bigr|
&\leqslant \max\!\Bigl\{\sup_{0\leqslant t\leqslant t_0}\bigl|\psi(\T_t(a)) - \varphi(\T_t(a))\bigr|,
\\[-1ex]
&\hspace{8em}\sup_{t\geqslant t_0}\bigl|\psi(\T_t(a)) - \varphi(\T_t(a))\bigr|\Bigr\}
<\varepsilon.
\end{align*}
Hence $(\T_t)$ has the $e$-property at $\varphi$. Since $\varphi$ was arbitrary, the theorem follows.
\end{proof}

\begin{remark}
The proof uses only norm continuity on $\A$. In particular, no asymptotic stability assumption is needed for Theorem~\ref{thm:eventual-implies-e}. It also does not require $\Ldom$ to be invariant under $(\T_t)$: for a fixed observable $a\in\Ldom$ the compactness argument only needs the orbit segment $\{\T_t(a):0\leqslant t\leqslant t_0\}$ to lie in the ambient Feller algebra $\A$.
\end{remark}

\begin{corollary}\label{prop:orbit-compactness}
Let $(\T_t)_{t\geqslant 0}$ be a quantum Markov semigroup on $\M$, and let $\A\subseteq\M$ be a Feller algebra. Fix $a\in\A_{\mathrm{sa}}$. Then the following are equivalent:
\begin{enumerate}[label=\upshape(\roman*)]
    \item the family $\{\T_t(a):t\geqslant 0\}$ is relatively compact in the norm topology of $\A$;
    \item the family of functions
    \[
    \Sst(\A)\ni\varphi\longmapsto \varphi(\T_t(a))\in \R,\qquad t\geqslant 0,
    \]
    is equicontinuous on the weak-$^*$ compact state space $\Sst(\A)$.
\end{enumerate}
Consequently, $(\T_t)$ has the $e$-property on a test space $\Ldom$ if and only if every $a\in\Ldom$ has norm relatively compact orbit in $\A$.
\end{corollary}

\begin{proof}
Let
\[
\Gamma\colon \A_{\mathrm{sa}}\longrightarrow C(\Sst(\A),\R),
\qquad
\Gamma(b)(\varphi):=\widehat b(\varphi):=\varphi(b).
\]
For self-adjoint $b$ one has
\[
\|b\|=\sup_{\varphi\in\Sst(\A)}|\varphi(b)|,
\]
so $\Gamma$ is an isometry. If the orbit $\{\T_t(a):t\geqslant 0\}$ is norm relatively compact, then Lemma~\ref{lem:compact-uniform} applied to its norm closure gives weak-$^*$ equicontinuity of the corresponding expectation functions. This proves \upshape(i)$\Rightarrow$\upshape(ii).

Conversely, assume \upshape(ii). The functions $\Gamma(\T_t(a))$ are uniformly bounded by $\|a\|$, because the maps $\T_t$ are unital and positive, hence contractive. Since $\Sst(\A)$ is weak-$^*$ compact, the Arzel\`a--Ascoli theorem implies that
\[
\{\Gamma(\T_t(a)):t\geqslant 0\}
\]
is relatively compact in $C(\Sst(\A),\R)$. Since $\Gamma$ is an isometry from the Banach space $\A_{\mathrm{sa}}$ into $C(\Sst(\A),\R)$, its range is closed. Thus relative compactness of the image implies relative compactness of the original orbit in $\A$.
\end{proof}

\begin{corollary}\label{cor:dense-test-extension}
Let $\Ldom\subseteq\A_{\mathrm{sa}}$ be norm dense. If $(\T_t)$ has the $e$-property on $\Ldom$, then it has the $e$-property on all of $\A_{\mathrm{sa}}$.
\end{corollary}

\begin{proof}
Fix $b\in\A_{\mathrm{sa}}$, $\varphi\in\Sst(\A)$, and $\varepsilon>0$. Choose $a\in\Ldom$ with $\|b-a\|<\varepsilon/6$. By the $e$-property for $a$, there exists a weak-$^*$ neighbourhood $U$ of $\varphi$ such that
\[
\sup_{\psi\in U}\sup_{t\geqslant 0}
|\psi(\T_t(a))-\varphi(\T_t(a))|<\varepsilon/3.
\]
For $\psi\in U$ and $t\geqslant 0$, contractivity of $\T_t$ gives
\begin{align*}
|\psi(\T_t(b))-\varphi(\T_t(b))|
&\leqslant |\psi(\T_t(a))-\varphi(\T_t(a))|
   +2\|\T_t(b-a)\| \\
&<\varepsilon/3+2\|b-a\|<\varepsilon.
\end{align*}
Thus the $e$-property holds for $b$.
\end{proof}

\begin{remark}
Corollary~\ref{prop:orbit-compactness} is specific to the full state-space formulation. In the classical Polish-space definition one usually tests equicontinuity only on the embedded Dirac states, which need not form a compact space. Passing to the whole compact state space converts the $e$-property into a norm-compactness statement for observable orbits. This is useful structurally, but it does not make the property vacuous: norm relative compactness of all relevant observable orbits is a strong dynamical restriction in infinite-dimensional algebras.
\end{remark}

\section{A quantitative full-state-space \texorpdfstring{$e$}{e}-property}\label{sec:quantitative}

\begin{definition}
Let $X$ be a compact Hausdorff space and let $\mathcal F\subseteq C(X,\R)$ be bounded.  Its \emph{equicontinuity defect} is
\[
   \eDef_X(\mathcal F)
   :=\inf_{\mathcal U}
      \max_{U\in\mathcal U}
      \sup_{f\in\mathcal F}\osc_U(f),
\]
where $\mathcal U$ ranges over finite open covers of $X$.  Its \emph{Hausdorff measure of non-compactness} is
\[
   \chi(\mathcal F)
   :=\inf\bigl\{r>0:\mathcal F\text{ has a finite }r\text{-net in }C(X,\R)\bigr\}.
\]
\end{definition}

The defect vanishes exactly for equicontinuous families.  The measure of non-compactness vanishes exactly for relatively compact families.

\begin{proposition}\label{prop:defect-noncompactness}
For every bounded family $\mathcal F\subseteq C(X,\R)$,
\begin{equation}\label{eq:defect-comparison}
   \frac12\eDef_X(\mathcal F)
   \leqslant\chi(\mathcal F)
   \leqslant\eDef_X(\mathcal F).
\end{equation}
If $\mathcal F$ and $\mathcal G$ are bounded, then
\begin{equation}\label{eq:defect-union}
   \eDef_X(\mathcal F\cup\mathcal G)
   =\max\{\eDef_X(\mathcal F),\eDef_X(\mathcal G)\}.
\end{equation}
\end{proposition}

\begin{proof}
Fix $r>\chi(\mathcal F)$ and choose $g_1,\ldots,g_m\in C(X,\R)$ such that every $f\in\mathcal F$ satisfies $\|f-g_j\|<r$ for some $j$.  The finite family $(g_j)$ is equicontinuous.  Hence, for every $\varepsilon>0$, there is a finite open cover $\mathcal U$ such that
\[
   \osc_U(g_j)<\varepsilon
   \qquad(U\in\mathcal U,\ 1\leqslant j\leqslant m).
\]
If $\|f-g_j\|<r$, then
\[
   \osc_U(f)\leqslant\osc_U(g_j)+2r<\varepsilon+2r.
\]
Thus $\eDef_X(\mathcal F)\leqslant2r$.  Letting $r\downarrow\chi(\mathcal F)$ proves the left inequality in \eqref{eq:defect-comparison}.

For the reverse inequality, fix $d>\eDef_X(\mathcal F)$.  Choose a finite open cover $\mathcal U=(U_j)_{j=1}^m$ such that
\[
   \sup_{f\in\mathcal F}\osc_{U_j}(f)<d
   \qquad(1\leqslant j\leqslant m).
\]
Choose a partition of unity $(h_j)_{j=1}^m$ subordinate to this cover and points $x_j\in U_j$ whenever $h_j\neq0$.  Define a finite-rank operator
\[
   Qf:=\sum_{j=1}^m f(x_j)h_j.
\]
If $h_j(x)>0$, then $x,x_j\in U_j$, and consequently
\[
   |f(x)-Qf(x)|
   \leqslant\sum_jh_j(x)|f(x)-f(x_j)|<d.
\]
Hence $\|f-Qf\|\leqslant d$ for $f\in\mathcal F$.  The set $Q(\mathcal F)$ is bounded in a finite-dimensional space, so it is relatively compact.  It has a finite $\varepsilon$-net for every $\varepsilon>0$, and therefore $\mathcal F$ has a finite $(d+\varepsilon)$-net.  Letting first $\varepsilon\downarrow0$ and then $d\downarrow\eDef_X(\mathcal F)$ proves the right inequality.

The lower bound in \eqref{eq:defect-union} is immediate.  For the converse, choose finite covers that realise the two defects up to $\varepsilon$ and take their common refinement.  On every member of the refinement, the oscillation of the union is bounded by the larger of the two estimates.  Let $\varepsilon\downarrow0$.
\end{proof}

\begin{remark}\label{rem:nussbaum-sharp}
Quantitative measures of non-equicontinuity go back to Nussbaum and have been developed in much greater generality; see \cite{Nussbaum1971,CaponettiTrombettaTrombetta2023}.  Proposition~\ref{prop:defect-noncompactness} is the scalar compact-Hausdorff form needed here.  Its factor two is sharp for general bounded families.  Indeed, if $X$ is infinite, the closed unit ball $B$ of $C(X,\R)$ satisfies $\chi(B)=1$ by Riesz's lemma.  Every finite open cover of $X$ contains a member with two distinct points, and normality of $X$ gives a function in $B$ taking the values $1$ and $-1$ at those points.  Hence $\eDef_X(B)=2$.
\end{remark}

Put
\[
   \mathcal O_{\T}(a):=\{\T_t(a):t\geqslant0\},
\]
and define
\[
   \chiT_{\T}(a):=\chi(\mathcal O_{\T}(a)),
   \qquad
   \eDef_{\T}(a):=
   \eDef_{\Sst(\A)}
   \bigl(\{\widehat{\T_t(a)}:t\geqslant0\}\bigr),
\]
where $\widehat b=\Gamma(b)$ denotes the evaluation function associated with $b\in\A_{\mathrm{sa}}$, and the first $\chi$ is computed in $\A_{\mathrm{sa}}$.

\begin{theorem}\label{thm:quantitative-e}
For $a,b\in\A_{\mathrm{sa}}$,
\begin{equation}\label{eq:q-e-comparison}
   \frac12\eDef_{\T}(a)
   \leqslant\chiT_{\T}(a)
   \leqslant\eDef_{\T}(a),
\end{equation}
\begin{equation}\label{eq:q-e-lipschitz}
\begin{split}
   |\chiT_{\T}(a)-\chiT_{\T}(b)|&\leqslant\|a-b\|,\\
   |\eDef_{\T}(a)-\eDef_{\T}(b)|&\leqslant2\|a-b\|.
\end{split}
\end{equation}
For every $T\geqslant0$, deleting the interval $[0,T]$ leaves both quantities unchanged.
\end{theorem}

\begin{proof}
Kadison's map
\[
   \Gamma\colon\A_{\mathrm{sa}}\longrightarrow C(\Sst(\A),\R),
   \qquad
   \Gamma(c)(\varphi)=\varphi(c),
\]
is an isometry.  Hence the Hausdorff measure of non-compactness of $\mathcal O_{\T}(a)$ equals that of its image under $\Gamma$.  Proposition~\ref{prop:defect-noncompactness} gives \eqref{eq:q-e-comparison}.

For every $t\geqslant0$,
\[
   \|\T_t(a)-\T_t(b)\|\leqslant\|a-b\|.
\]
Thus the two orbit sets are at Hausdorff distance at most $\|a-b\|$, which proves the first estimate in \eqref{eq:q-e-lipschitz}.  On every subset $U\subseteq\Sst(\A)$,
\[
   \osc_U(\widehat{\T_t(a)})
   \leqslant\osc_U(\widehat{\T_t(b)})+2\|a-b\|.
\]
Taking the relevant suprema and infima, and then interchanging $a$ and $b$, gives the second estimate.

Strong continuity makes the initial orbit segment
\[
   K_T:=\{\T_t(a):0\leqslant t\leqslant T\}
\]
norm compact.  The Hausdorff measure of non-compactness has the maximum property under finite unions, so $\chi(K_T)=0$ and adjoining $K_T$ to the tail does not change its value.  The family $\Gamma(K_T)$ is equicontinuous, so its defect is zero.  Formula \eqref{eq:defect-union} shows that adjoining it to the tail does not change the equicontinuity defect.
\end{proof}

\begin{corollary}\label{cor:e-equivalences}
For $a=a^*\in\A$, the following are equivalent:
\begin{enumerate}[label=\upshape(\roman*)]
  \item the full-state-space $e$-property holds for $a$;
  \item the eventual full-state-space $e$-property holds for $a$;
  \item $\mathcal O_{\T}(a)$ is relatively compact in norm;
  \item $\eDef_{\T}(a)=0$;
  \item $\chiT_{\T}(a)=0$.
\end{enumerate}
\end{corollary}

\begin{proof}
The equivalence of \upshape(i), \upshape(iv) and \upshape(v) follows from Theorem~\ref{thm:quantitative-e}.  The full property implies the eventual one.  Conversely, assume the eventual property.  Fix $\varepsilon>0$.  For every state $\varphi$, choose a neighbourhood $U_\varphi$ and a time $T_\varphi$ giving tail oscillation less than $\varepsilon$.  Compactness of the state space gives a finite subcover.  Taking the maximum of the corresponding times shows that a single tail has equicontinuity defect at most $2\varepsilon$; the harmless factor two comes from comparing two states in one neighbourhood with its centre.  By Theorem~\ref{thm:quantitative-e}, deleting the initial segment does not change the defect of the full family.  Since $\varepsilon$ is arbitrary, the full defect is zero.
\end{proof}

\begin{remark}
On a compact classical phase space $X$, Corollary~\ref{cor:e-equivalences} recovers the ordinary $e$-property because Dirac states identify $X$ with a compact subset of $\Sst(C(X))$ and classical Arzel\`a--Ascoli gives the converse.  On a non-compact space, testing the full state space of $C_b(X)$ is strictly stronger than testing only the Dirac states.
\end{remark}

\begin{example}\label{ex:translation-failure}
Let $\A=\operatorname{BUC}(\R)$ and
\[
   (\T_t f)(x)=f(x+t).
\]
This is a strongly continuous group of unital $*$-automorphisms.  For $h(x)=e^{-x^2}$, the sequence $(\T_nh)_{n\geqslant1}$ is separated in norm: if $m-n$ is large, evaluation at $-n$ gives
\[
   \|\T_nh-\T_mh\|\geqslant1-e^{-(m-n)^2}.
\]
Thus the orbit is not relatively compact and the full-state-space $e$-property fails for $h$.
\end{example}

\subsection{The hierarchy of eventual equicontinuity variants}\label{subsec:hierarchy}

To place Corollary~\ref{cor:e-equivalences} in the context of the classical comparison theory \cite{Czapla2018,HilleSzarekWormZiemlanska2021,KukulskiWojewodka2024,LiuLiu2024}, we record the natural weak-$^*$ variants of eventual equicontinuity at a state $\varphi\in\Sst(\A)$, where $\mathcal U(\varphi)$ denotes the family of weak-$^*$ neighbourhoods of $\varphi$ in $\Sst(\A)$.

\begin{definition}\label{def:equicontinuity-variants}
Let $a\in\A_{\mathrm{sa}}$ and $\varphi\in\Sst(\A)$.
\begin{enumerate}[label=\upshape(\arabic*)]
    \item $(\T_t)$ has the \emph{fixed-threshold eventual $e$-property} at $\varphi$ for $a$ if there exists $t_0\geqslant 0$ such that for every $\varepsilon>0$ there exists $U\in\mathcal U(\varphi)$ satisfying
    \[
    \sup_{\psi\in U}\sup_{t\geqslant t_0}\bigl|\psi(\T_t(a)) - \varphi(\T_t(a))\bigr| < \varepsilon.
    \]
    \item $(\T_t)$ is \emph{locally uniformly eventually continuous} at $\varphi$ for $a$ if
    \[
    \inf_{U\in\mathcal U(\varphi)}\limsup_{t\to\infty}\sup_{\psi\in U}\bigl|\psi(\T_t(a)) - \varphi(\T_t(a))\bigr| = 0.
    \]
    \item $(\T_t)$ is \emph{pointwise eventually continuous} at $\varphi$ for $a$ if
    \[
    \inf_{U\in\mathcal U(\varphi)}\sup_{\psi\in U}\limsup_{t\to\infty}\bigl|\psi(\T_t(a)) - \varphi(\T_t(a))\bigr| = 0.
    \]
\end{enumerate}
If the corresponding condition holds at every $\varphi\in\Sst(\A)$, we say that it holds on $\Sst(\A)$ for $a$.
\end{definition}

\begin{remark}\label{rem:quantifier-order}
The distinction between \upshape(2) and \upshape(3) lies entirely in the order of $\limsup_{t\to\infty}$ and $\sup_{\psi\in U}$. Condition~\upshape(2) requires a single tail threshold to control an entire weak-$^*$ neighbourhood $U$ up to tolerance $\varepsilon$, which is filter language for the $\varepsilon$-dependent eventual $e$-property of Section~\ref{sec:framework}. Condition~\upshape(3) is the direct state-space analogue of the pointwise classical condition $\lim_{x\to z}\limsup_{t\to\infty}|P_t f(x)-P_t f(z)|=0$ studied in \cite{LiuLiu2024}. As we show below, Condition~\upshape(2) collapses completely to the ordinary $e$-property on $\Sst(\A)$, whereas Condition~\upshape(3) is strictly weaker.
\end{remark}

\begin{theorem}\label{thm:collapse}
Let $(\T_t)_{t\geqslant 0}$ be a quantum Markov semigroup on $\M$, let $\A\subseteq\M$ be a Feller algebra, and fix $a\in\A_{\mathrm{sa}}$. The following assertions are equivalent:
\begin{enumerate}[label=\upshape(\roman*)]
    \item $(\T_t)$ has the ordinary $e$-property on $\Sst(\A)$ for $a$;
    \item $(\T_t)$ has the fixed-threshold eventual $e$-property on $\Sst(\A)$ for $a$;
    \item $(\T_t)$ has the eventual $e$-property on $\Sst(\A)$ for $a$;
    \item $(\T_t)$ is locally uniformly eventually continuous on $\Sst(\A)$ for $a$;
    \item the forward orbit $\mathcal O_\T(a)=\{\T_t(a):t\geqslant 0\}$ is norm relatively compact in $\A$;
    \item $\eDef_\T(a)=0$;
    \item $\chiT_\T(a)=0$.
\end{enumerate}
Each of these conditions implies that $(\T_t)$ is pointwise eventually continuous on $\Sst(\A)$ for $a$, but the converse fails in general.
\end{theorem}

\begin{proof}
The equivalences $(i)\Longleftrightarrow(iii)\Longleftrightarrow(v)\Longleftrightarrow(vi)\Longleftrightarrow(vii)$ are given by Corollary~\ref{cor:e-equivalences}, and $(i)\Longrightarrow(ii)\Longrightarrow(iii)$ is immediate.

$(iii)\Longrightarrow(iv)$:  Fix $\varphi\in\Sst(\A)$ and $\varepsilon>0$. By $(iii)$, there exist $t_0\geqslant 0$ and $U\in\mathcal U(\varphi)$ such that
\[
\sup_{\psi\in U}\sup_{t\geqslant t_0}\bigl|\psi(\T_t(a))-\varphi(\T_t(a))\bigr|\leqslant \varepsilon.
\]
Hence
\[
\limsup_{t\to\infty}\sup_{\psi\in U}\bigl|\psi(\T_t(a))-\varphi(\T_t(a))\bigr|\leqslant \varepsilon.
\]
Taking the infimum over all $U\in\mathcal{U}(\varphi)$ yields an infimum bounded by $\varepsilon$. Since $\varepsilon>0$ was arbitrary, $\mathrm{(iv)}$ holds.

$(iv)\Longrightarrow(vi)$: Let $f_t:=\Gamma(\T_t(a))\in C(\Sst(\A),\R)$, where $\Gamma$ is the evaluation map used in Corollary~\ref{prop:orbit-compactness}, and fix $\varepsilon>0$. By $(iv)$, for every $\varphi\in\Sst(\A)$ there exist an open weak-$^*$ neighbourhood $U_\varphi$ of $\varphi$ and a time $\tau_\varphi\geqslant 0$ such that
\[
\sup_{\psi\in U_\varphi}\bigl|f_t(\psi)-f_t(\varphi)\bigr|\leqslant \frac{\varepsilon}{2}
\qquad (t\geqslant \tau_\varphi).
\]

By the triangle inequality, $\osc_{U_\varphi}(f_t)\leqslant \varepsilon$ for all $t\geqslant \tau_\varphi$. By weak-$^*$ compactness of $\Sst(\A)$, choose a finite subcover $\{U_{\varphi_1},\dots,U_{\varphi_m}\}$ and put $\tau_{\max}:=\max_{1\leqslant j\leqslant m}\tau_{\varphi_j}$. Then the tail family $\mathcal F_{\tau_{\max}}:=\{f_t:t\geqslant \tau_{\max}\}$ satisfies
\[
\eDef_{\Sst(\A)}(\mathcal F_{\tau_{\max}})
\leqslant \max_{1\leqslant j\leqslant m}\sup_{t\geqslant \tau_{\max}}\osc_{U_{\varphi_j}}(f_t)
\leqslant \varepsilon.
\]
By Theorem~\ref{thm:quantitative-e}, deleting the initial interval $[0,\tau_{\max}]$ leaves the equicontinuity defect unchanged, so $\eDef_\T(a)=\eDef_{\Sst(\A)}(\mathcal F_{\tau_{\max}})\leqslant\varepsilon$. Since $\varepsilon>0$ was arbitrary, $\eDef_\T(a)=0$.

\end{proof}

\section{Ces\`aro averages and the mean-ergodic mechanism}\label{sec:cesaro}

For $a\in\A$ and $T>0$ we put
\[
M_T(a):=\frac1T\int_0^T \T_t(a)\,dt,
\]
where the integral is the Bochner integral in $\A$.

\begin{definition}
Let $(\T_t)_{t\geqslant 0}$ be a quantum Markov semigroup on $\M$ and let $\A\subseteq\M$ be a Feller algebra. We say that $(\T_t)$ has the \emph{Ces\`aro $e$-property at $\varphi\in\Sst(\A)$ relative to $\Ldom$} if for every $a\in\Ldom$ and every $\varepsilon>0$ there exists a weak-$^*$ neighbourhood $U$ of $\varphi$ such that
\[
\sup_{\psi\in U}\sup_{T>0}\bigl|\psi(M_T(a)) - \varphi(M_T(a))\bigr|<\varepsilon.
\]
If this holds at every state, we say that $(\T_t)$ has the Ces\`aro $e$-property on $\Ldom$.
\end{definition}

\begin{lemma}\label{lem:mean-cont}
For every $a\in\A$, the map $(0,\infty)\ni T\mapsto M_T(a)$ extends continuously to $T=0$ by the value $M_0(a):=a$.
\end{lemma}

\begin{proof}
For $T>0$ the map is norm continuous by standard properties of the Bochner integral and the norm continuity of $t\mapsto \T_t(a)$. Moreover,
\[
\|M_T(a)-a\|\leqslant \frac1T\int_0^T \|\T_t(a)-a\|\,dt\longrightarrow 0
\qquad (T\downarrow 0),
\]
again by norm continuity.
\end{proof}

\begin{proposition}\label{prop:e-implies-cesaro}
If $(\T_t)$ has the $e$-property on $\Ldom$, then it has the Ces\`aro $e$-property on $\Ldom$.
\end{proposition}

\begin{proof}
Fix $a\in\Ldom$, $\varphi\in\Sst(\A)$, and $\varepsilon>0$. By the $e$-property there exists a weak-$^*$ neighbourhood $U$ of $\varphi$ such that
\[
\sup_{\psi\in U}\sup_{t\geqslant 0}\bigl|\psi(\T_t(a)) - \varphi(\T_t(a))\bigr|<\varepsilon.
\]
Therefore, for every $\psi\in U$ and every $T>0$,
\begin{align*}
\bigl|\psi(M_T(a)) - \varphi(M_T(a))\bigr|
&\leqslant \frac1T\int_0^T \bigl|\psi(\T_t(a)) - \varphi(\T_t(a))\bigr|\,dt \\
&<\varepsilon.
\end{align*}
Hence the Ces\`aro $e$-property holds.
\end{proof}

\begin{theorem}\label{thm:mean-ergodic}
Let $(\T_t)_{t\geqslant 0}$ be a quantum Markov semigroup on $\M$, and let $\A\subseteq\M$ be a Feller algebra. Assume that for every $a\in\Ldom$ there exists $E(a)\in\A$ such that
\[
\|M_T(a)-E(a)\|\longrightarrow 0
\qquad (T\to\infty).
\]
Then for every $s\geqslant 0$ and every $a\in\Ldom$ one has
\[
\T_s(E(a))=E(a).
\]
In particular, $(\T_t)$ has the Ces\`aro $e$-property on $\Ldom$.
\end{theorem}

\begin{proof}
Fix $a\in\Ldom$ and $s\geqslant 0$. For $T>0$ we have
\[
\T_s(M_T(a))-M_T(a)
=
\frac1T\int_T^{T+s}\T_r(a)\,dr
-
\frac1T\int_0^s\T_r(a)\,dr,
\]
whence, using the contractivity of the positive unital maps $\T_r$,
\[
\|\T_s(M_T(a))-M_T(a)\|\leqslant \frac{2s\|a\|}{T}.
\]
Therefore
\begin{align*}
\|\T_s(E(a))-E(a)\|
&\leqslant \|\T_s(E(a)-M_T(a))\| + \|\T_s(M_T(a))-M_T(a)\| + \|M_T(a)-E(a)\| \\
&\leqslant 2\|M_T(a)-E(a)\| + \frac{2s\|a\|}{T}.
\end{align*}
Letting $T\to\infty$ gives $\T_s(E(a))=E(a)$.

Fix now $\varphi\in\Sst(\A)$ and $\varepsilon>0$. Choose $T_0>0$ such that
\[
\|M_T(a)-E(a)\|<\varepsilon/3,
\qquad T\geqslant T_0.
\]
Since evaluation at the fixed element $E(a)$ is weak-$^*$ continuous, there exists a weak-$^*$ neighbourhood $U_1$ of $\varphi$ such that
\[
|\psi(E(a)) - \varphi(E(a))|<\varepsilon/3,
\qquad \psi\in U_1.
\]
By Lemma~\ref{lem:mean-cont}, the set
\[
K:=\{M_T(a):0\leqslant T\leqslant T_0\}
\]
is norm compact. Lemma~\ref{lem:compact-uniform} therefore yields a weak-$^*$ neighbourhood $U_2$ of $\varphi$ such that
\[
\sup_{\psi\in U_2}\sup_{b\in K}|\psi(b)-\varphi(b)|<\varepsilon.
\]
Let $U:=U_1\cap U_2$. If $\psi\in U$ and $0<T\leqslant T_0$, then $M_T(a)\in K$, so
\[
|\psi(M_T(a)) - \varphi(M_T(a))|<\varepsilon.
\]
If instead $T\geqslant T_0$, then
\begin{align*}
|\psi(M_T(a)) - \varphi(M_T(a))|
&\leqslant |\psi(M_T(a)-E(a))| \\
&\quad + |\psi(E(a)) - \varphi(E(a))| \\
&\quad + |\varphi(E(a)-M_T(a))| \\
&\leqslant 2\|M_T(a)-E(a)\| + |\psi(E(a)) - \varphi(E(a))| \\
&<\varepsilon.
\end{align*}
Hence the Ces\`aro $e$-property holds at $\varphi$, and therefore on all of $\Sst(\A)$.
\end{proof}

\begin{corollary}\label{cor:scalar-mean-ergodic}
Under the assumptions of Theorem~\ref{thm:mean-ergodic}, assume in addition that there exists a state $\varphi_*\in\Sst(\A)$ such that
\[
\|M_T(a)-\varphi_*(a)1\|\longrightarrow 0
\qquad (T\to\infty)
\]
for every $a\in\Ldom$. Then for every $a\in\Ldom$,
\[
\sup_{\varphi,\psi\in\Sst(\A)}\bigl|\varphi(M_T(a)) - \psi(M_T(a))\bigr|\leqslant 2\|M_T(a)-\varphi_*(a)1\|\longrightarrow 0.
\]
In particular, the Ces\`aro $e$-property holds globally and uniformly.
\end{corollary}

\begin{proof}
For any states $\varphi,\psi$,
\begin{align*}
\bigl|\varphi(M_T(a)) - \psi(M_T(a))\bigr|
&\leqslant \bigl|\varphi(M_T(a)-\varphi_*(a)1)\bigr| + \bigl|\psi(M_T(a)-\varphi_*(a)1)\bigr| \\
&\leqslant 2\|M_T(a)-\varphi_*(a)1\|.
\end{align*}
The conclusion follows.
\end{proof}

\begin{remark}
Theorem~\ref{thm:mean-ergodic} is the abstract place where a $C^*$-algebraic von Neumann mean-ergodic theorem enters the picture. Whenever one has norm convergence of the ergodic averages to the mean-ergodic projection onto the fixed-point space, one obtains the Ces\`aro $e$-property automatically. Under the usual additional hypotheses, this fixed-point space may be identified with a fixed-point algebra. This is weaker than the ordinary $e$-property, but it is often easier to verify and it is already sufficient for uniform control of time-averaged expectations.
\end{remark}

\subsection{Quantitative Ces\`aro defects and strict separation from the \texorpdfstring{$e$}{e}-property}\label{subsec:cesaro-separation}

For $a\in\A_{\mathrm{sa}}$, write $\mathcal O_\T^{\mathrm{Ces}}(a):=\{M_T(a):T>0\}$ and define the \emph{Ces\`aro non-compactness measure} and \emph{Ces\`aro equicontinuity defect} by
\[
\chiT_\T^{\mathrm{Ces}}(a):=\chi\bigl(\mathcal O_\T^{\mathrm{Ces}}(a)\bigr),
\qquad
\eDef_\T^{\mathrm{Ces}}(a):=\eDef_{\Sst(\A)}\bigl(\{\widehat{M_T(a)}:T>0\}\bigr).
\]

\begin{proposition}\label{prop:quantitative-cesaro}
For every $a\in\A_{\mathrm{sa}}$,
\[
\frac12\eDef_\T^{\mathrm{Ces}}(a)\leqslant \chiT_\T^{\mathrm{Ces}}(a)\leqslant \eDef_\T^{\mathrm{Ces}}(a),
\]
and
\[
\chiT_\T^{\mathrm{Ces}}(a)\leqslant \chiT_\T(a),
\qquad
\eDef_\T^{\mathrm{Ces}}(a)\leqslant \eDef_\T(a).
\]
\end{proposition}

\begin{proof}
The first chain of inequalities follows from Proposition~\ref{prop:defect-noncompactness} applied to the family $\{\widehat{M_T(a)}:T>0\}$. For every subset $U\subseteq\Sst(\A)$ and $T>0$,
\[
\osc_U(\widehat{M_T(a)})
= \sup_{\varphi,\psi\in U}\left|\frac1T\int_0^T\bigl(\varphi(\T_t(a))-\psi(\T_t(a))\bigr)\,dt\right|
\leqslant \sup_{t\geqslant 0}\osc_U(\widehat{\T_t(a)}),
\]
which gives $\eDef_\T^{\mathrm{Ces}}(a)\leqslant \eDef_\T(a)$.

To compare the measures of non-compactness, fix $r>\chiT_\T(a)$ and choose a finite $r$-net $F=\{b_1,\dots,b_m\}\subseteq\A_{\mathrm{sa}}$ for $\mathcal O_\T(a)$. Then $\mathcal O_\T(a)$ is contained in the closed convex set $\operatorname{conv}(F)+\overline{B}_{\A_{\mathrm{sa}}}(0,r)$. Since $M_T(a)$ is a limit of Bochner averages of elements of $\mathcal O_\T(a)$, we have
\[
\mathcal O_\T^{\mathrm{Ces}}(a)\subseteq \operatorname{conv}(F)+\overline{B}_{\A_{\mathrm{sa}}}(0,r).
\]
Because $\operatorname{conv}(F)$ is compact in a finite-dimensional subspace, it admits a finite $\varepsilon$-net for every $\varepsilon>0$, and hence $\mathcal O_\T^{\mathrm{Ces}}(a)$ admits a finite $(r+\varepsilon)$-net. Letting $\varepsilon\downarrow 0$ and $r\downarrow\chiT_\T(a)$ yields $\chiT_\T^{\mathrm{Ces}}(a)\leqslant\chiT_\T(a)$.
\end{proof}

We now show that both inequalities in Proposition~\ref{prop:quantitative-cesaro} can be strict ($\chiT_\T^{\mathrm{Ces}}(a)=\eDef_\T^{\mathrm{Ces}}(a)=0$ while $\eDef_\T(a)>0$), both for discrete-time matrix Jamison systems and for continuous-time quantum Markov--Feller semigroups.

\begin{proposition}[Matrix Jamison System]\label{prop:jamison}
Let $S=\{0\}\cup\{1/n:n\in\N\}\subseteq [0,1]$ with the Euclidean metric, and let $\A=C(S,M_m(\C))$ with $m\geqslant 1$. We identify each $b\in M_m(\C)$ with the constant function $x\mapsto b$ in $\A$. Define $T\colon S\to S$ by
\[
T(0)=0,\qquad T(1)=0,\qquad T(1/n)=\frac{1}{n-1}\quad(n\geqslant 2),
\]
and let $P\colon \A\to\A$ be the unital completely positive map $P(f)=f\circ T$. Then:
\begin{enumerate}[label=\upshape(\roman*)]
    \item the discrete ergodic averages $M_N(P)(f) = \frac1N\sum_{k=0}^{N-1}P^k(f)$ satisfy
    \[
    \lim_{N\to\infty} \bigl\|M_N(P)(f) - f(0)\bigr\|_\A = 0 \qquad(f\in\A),
    \]
    so $P$ satisfies the discrete Ces\`aro $e$-property.
    \item the discrete semigroup $(P^k)_{k\geqslant 0}$ fails the ordinary $e$-property on $\A_{\mathrm{sa}}$.
\end{enumerate}
\end{proposition}

\begin{proof}
\upshape(i) Let $f\in\A$ with $f(0)=0$. Fix $\varepsilon>0$ and choose $n_0\in\N$ such that $\|f(1/n)\|_{M_m}<\varepsilon/2$ for all $n\geqslant n_0$. For any $n,N\in\N$, the iterates satisfy $T^k(1/n)=1/(n-k)$ for $k<n$ and $T^k(1/n)=0$ for $k\geqslant n$. Since $f(0)=0$, we have $f(T^k(1/n))=0$ for all $k\geqslant n$. Therefore,
\[
\left\|\sum_{k=0}^{N-1} f(T^k(1/n))\right\|_{M_m}
= \left\|\sum_{j=\max\{1,n-N+1\}}^n f(1/j)\right\|_{M_m}
\leqslant n_0\|f\|_\infty + N\frac{\varepsilon}{2},
\]
because at most $n_0$ indices satisfy $j\leqslant n_0$. Dividing by $N$, we obtain for all $N\geqslant 2n_0\|f\|_\infty/\varepsilon$ that
\[
\sup_{x\in S} \left\|\frac1N\sum_{k=0}^{N-1} f(T^k(x))\right\|_{M_m}
\leqslant \frac{n_0\|f\|_\infty}{N} + \frac{\varepsilon}{2} \leqslant \varepsilon.
\]
Thus $\|M_N(P)(f)\|_\A \to 0$. For general $f\in\A$, decomposing $f = (f-f(0)) + f(0)$ gives norm convergence of $M_N(P)(f)$ to the constant function $f(0)$, which yields the discrete Ces\`aro $e$-property.

\upshape(ii) Choose $g\in C(S,\R)$ such that $g(0)=0$, $g(1)=1$, and $0\leqslant g\leqslant 1$. Let $e_{11}\in M_m(\C)$ be a rank-one projection and set $a(x)=g(x)e_{11}\in\A_{\mathrm{sa}}$. Define states $\varphi_n(f)=\operatorname{tr}_m(f(1/n))$ and $\varphi_0(f)=\operatorname{tr}_m(f(0))$, where $\operatorname{tr}_m$ denotes the normalised trace. Then $\varphi_n\to\varphi_0$ weak-$^*$ in $\Sst(\A)$ as $n\to\infty$. For each $n\geqslant 2$, taking $k=n-1$ gives $T^{n-1}(1/n)=1$, so
\[
\varphi_n(P^{n-1}(a)) = \operatorname{tr}_m(a(1)) = \operatorname{tr}_m(e_{11}) = \frac{1}{m},
\]
whereas $\varphi_0(P^k(a)) = \operatorname{tr}_m(a(0)) = 0$ for all $k\geqslant 0$. Thus
\[
\sup_{k\geqslant 0}\bigl|\varphi_n(P^k(a)) - \varphi_0(P^k(a))\bigr| \geqslant \frac{1}{m} > 0
\qquad (n\geqslant 2),
\]
and the ordinary $e$-property fails at $\varphi_0$.
\end{proof}

\begin{example}[Continuous-time homoclinic matrix flow]\label{ex:continuous-homoclinic}
Identify the unit circle $\mathbb T$ with $\R/2\pi\mathbb Z$. Let $\mu$ be normalised Haar measure on $\mathbb T$, let $\mathfrak B(\mathbb T)$ be its Borel $\sigma$-algebra, and fix an integer $m\geqslant1$. Set
\[
\M = L^\infty(\mathbb T,\mathfrak B(\mathbb T),\mu)\,\bar\otimes\,M_m(\C),
\qquad
\A = C(\mathbb T,M_m(\C)).
\]
Let $(\phi_t)_{t\in\R}$ be the smooth flow on $\mathbb T$ generated by the vector field $\dot\theta = 1-\cos\theta = 2\sin^2(\theta/2)\geqslant 0$, and define
\[
(\T_t f)(z) := f(\phi_t(z))
\qquad (f\in\M,\ t\geqslant 0).
\]
Because each $\phi_t$ is a diffeomorphism of $\mathbb T$ preserving the measure class of $\mu$, $(\T_t)_{t\geqslant 0}$ is a weak-$^*$ continuous normal quantum Markov semigroup on $\M$ for which $\A$ is a Feller algebra.
\begin{enumerate}[label=\upshape(\roman*)]
    \item The point $z=1$ ($\theta=0$) is a global attractor: $\lim_{t\to\infty}\phi_t(z)=1$ for every $z\in\mathbb T$. Given $f\in\A$ with $f(1)=0$ and $\varepsilon>0$, choose an open arc $U\subseteq\mathbb T$ around $1$ on which $\|f(z)\|_{M_m}<\varepsilon/2$. On the compact complement $\mathbb T\setminus U$, the angular velocity $2\sin^2(\theta/2)$ is bounded below by a constant $c_U>0$, so every forward trajectory spends total time at most $\tau_U:=2\pi/c_U<\infty$ in $\mathbb T\setminus U$. Consequently, for every $T>0$,
\begin{align*}
\|M_T(f)\|_\A
&\leqslant \sup_{z\in\mathbb T}\left(
   \frac{1}{T}\int_{I_{\mathrm{out}}(z)} \bigl\Vert{}f(\phi_t(z))\bigr\Vert{}_{M_m}\,dt
   + \frac{1}{T}\int_{I_{\mathrm{in}}(z)} \bigl\Vert{}f(\phi_t(z))\bigr\Vert{}_{M_m}\,dt
\right) \\
&\leqslant \frac{\tau_U\|f\|_\infty}{T} + \frac{\varepsilon}{2}
,
\end{align*}
    where $I_{\mathrm{in}}(z) := \{t \in [0, T] : \phi_t(z) \in U\}$ and $I_{\mathrm{out}}(z) := \{t \in [0, T] : \phi_t(z) \in \mathbb{T} \setminus U\}$. 
    Taking the limit superior as $T\to\infty$ and then letting $\varepsilon\downarrow0$ gives $\|M_T(f)\|_\A\to0$. Identifying matrices with constant functions, we obtain $M_T(f)\to f(1)$ in norm for every $f\in\A$. By Theorem~\ref{thm:mean-ergodic}, $(\T_t)$ satisfies the continuous-time Ces\`aro $e$-property on $\A_{\mathrm{sa}}$ (with $\chiT_\T^{\mathrm{Ces}}(a)=\eDef_\T^{\mathrm{Ces}}(a)=0$ for all $a\in\A_{\mathrm{sa}}$).
    \item For $a(e^{i\theta}):=\sin^2(\theta/2)\,e_{11}\in\A_{\mathrm{sa}}$, the states $\varphi_n(f):=\operatorname{tr}_m(f(e^{i/n}))$ converge weak-$^*$ to $\varphi_0(f):=\operatorname{tr}_m(f(1))$, yet for each $n\in\N$ there exists $t_n>0$ with $\phi_{t_n}(e^{i/n})=-1$, giving
    \[
    \sup_{t\geqslant 0}\bigl|\varphi_n(\T_t(a))-\varphi_0(\T_t(a))\bigr| \geqslant \varphi_n(\T_{t_n}(a)) = \frac1m > 0.
    \]
    Hence $(\T_t)$ fails the ordinary $e$-property on $\A_{\mathrm{sa}}$.
\end{enumerate}
\end{example}

\section{Reversible and stable dynamics}\label{sec:reversible}

Write
\[
   \AP(\T):=\{a\in\A:\mathcal O_{\T}(a)\text{ is relatively compact in norm}\}.
\]
It is a closed invariant operator system.  If every $\T_t$ is a $*$-endomorphism, it is an ambient $C^*$-subalgebra.  For general completely positive dynamics, the reversible range carries a different product.

The point-norm topology on a family of maps is the topology of pointwise norm convergence; we abbreviate it by $\pn$. We also write
\[
   \Fix(\T):=\{a\in\A:\T_t(a)=a\text{ for every }t\geqslant0\},
\]
and
\[
   \Inv(\T):=\{\varphi\in\Sst(\A):\varphi\circ\T_t=\varphi\text{ for every }t\geqslant0\}.
\]

\begin{theorem}\label{thm:jdlg}
Let $(\T_t)_{t\geqslant0}$ be a quantum Markov semigroup on $\M$, let $\A\subseteq\M$ be a Feller algebra, and suppose that $\AP(\T)=\A$.  Then there is a unital completely positive projection $P\colon\A\to\A$ with the following properties:
\begin{enumerate}[label=\upshape(\roman*)]
  \item $P\T_t=\T_tP$ for every $t\geqslant0$;
  \item
  \[
      \A=\Rr\oplus\ker P,
      \qquad \Rr=P(\A);
  \]
  \item $\|\T_t(a)\|\to0$ for every $a\in\ker P$;
  \item
  \[
      G:=\overline{\{\T_t|_{\Rr}:t\geqslant0\}}^{\pn}
  \]
  is a compact abelian group of surjective complete order isomorphisms of $\Rr$;
  \item with
  \[
      x\circ y:=P(xy),
      \qquad x,y\in\Rr,
  \]
  the operator system $\Rr$ is a $C^*$-algebra and every $g\in G$ is a $*$-automorphism for this product.
\end{enumerate}
\end{theorem}

\begin{proof}
By Corollary~\ref{cor:e-equivalences}, all Banach-space orbits are norm relatively compact.  The Jacobs--de Leeuw--Glicksberg theorem for strongly almost periodic semigroups gives a bounded projection $P$ in the point-norm closure of the semigroup, commuting with it, a reversible range and a stable kernel; see \cite{deLeeuwGlicksberg1961,Eisner2010}.  On the stable kernel, zero belongs to the norm closure of every orbit.  Contractivity upgrades this to convergence of the whole orbit: given $\varepsilon>0$, choose $s$ with $\|\T_s(a)\|<\varepsilon$; then
\[
   \|\T_t(a)\|
   =\|\T_{t-s}(\T_s(a))\|
   \leqslant\|\T_s(a)\|<\varepsilon
   \qquad(t\geqslant s).
\]
This proves \upshape(iii) explicitly.

Point-norm limits of unital completely positive maps are unital and completely positive, so $P$ has these properties.  The closure of a commuting semigroup is abelian.  Every element of $G$ and its inverse is a point-norm limit of restrictions of unital completely positive maps; hence the elements of $G$ are complete order automorphisms.  The Choi--Effros theorem gives the $C^*$-product $x\circ y=P(xy)$ on the range \cite{ChoiEffros1977}, and complete order automorphisms are $*$-automorphisms for that product.
\end{proof}

\begin{remark}
Theorem~\ref{thm:jdlg} is a $C^*$-algebraic, point-norm formulation of established Jacobs--de Leeuw--Glicksberg and decoherence theory.  Hellmich treated weak-$^*$ continuous quantum dynamical semigroups on dual spaces \cite{Hellmich2011}; B\'atkai, Groh, Kunszenti-Kov\'acs and Schreiber obtained structural decompositions on $W^*$-algebras \cite{BatkaiGrohKunszentiSchreiber2012}; Runde and Viselter developed almost-periodic operators for quantum semigroups \cite{RundeViselter2014}.  Earlier long-time asymptotics were established by Frigerio and Frigerio--Verri \cite{Frigerio1978,FrigerioVerri1982}, while modern decoherence-free decompositions include \cite{CarboneSassoUmanita2014,FagnolaSassoUmanita2019}.  The Choi--Effros and peripheral viewpoints are further developed in \cite{FidaleoOttomanoRossi2022,BhatKarTalwar2023,BhatDias2025,AmatoFacchiKonderak2026}.  We do not claim this architecture as new; it is the structural setting for the quantitative and multiplicative results of the present paper.
\end{remark}

\begin{theorem}\label{thm:invariant-states}
Let $P$, $\Rr$ and $G$ be as in Theorem~\ref{thm:jdlg}.  Restriction induces an affine weak-$^*$ homeomorphism
\[
   \Inv(\T)
   \longrightarrow
   \Sst(\Rr)^G
   :=\{\psi\in\Sst(\Rr):\psi\circ g=\psi\text{ for every }g\in G\}.
\]
Its inverse is $\psi\mapsto\psi\circ P$.
\end{theorem}

\begin{proof}
Let $\varphi$ be invariant.  For $a\in\ker P$,
\[
   \varphi(a)=\varphi(\T_t(a))\longrightarrow0,
\]
so $\varphi=\varphi\circ P$.  Its restriction to $\Rr$ is invariant under every $\T_t|_{\Rr}$ and hence, by point-norm continuity, under $G$.

Conversely, if $\psi\in\Sst(\Rr)^G$, then $\psi\circ P$ is a state on $\A$ because $P$ is unital and positive.  Since $P$ commutes with $\T_t$ and the restriction of $\T_t$ belongs to $G$,
\[
   (\psi\circ P)\circ\T_t
   =\psi\circ(\T_t|_{\Rr})\circ P
   =\psi\circ P.
\]
The two maps are inverse.  Weak-$^*$ continuity of restriction and composition with $P$ gives the homeomorphism.
\end{proof}

For $\lambda\in\R$, put
\[
   \A_\lambda
   :=\{a\in\A:\T_t(a)=e^{i\lambda t}a\text{ for every }t\geqslant0\}.
\]

\begin{theorem}\label{thm:spectral-synthesis}
Under the assumptions of Theorem~\ref{thm:jdlg},
\[
   \Rr=
   \overline{\operatorname{span}}
   \{\A_\lambda:\lambda\in\R\}.
\]
\end{theorem}

\begin{proof}
For a character $\gamma\in\widehat G$, define the Fourier projection
\[
   P_\gamma a
   :=\int_G\overline{\gamma(g)}\,g(a)\,dm_G(g),
   \qquad a\in\Rr,
\]
where $m_G$ is normalised Haar measure.  A net of trigonometric approximate identities on the compact abelian group $G$ converges uniformly on every continuous orbit map $g\mapsto g(a)$; see, for example, \cite[Chapter~2]{Rudin1962}.  Consequently, the linear span of the Fourier subspaces $P_\gamma(\Rr)$ is dense in $\Rr$.

The one-parameter group $t\mapsto\T_t|_{\Rr}$ has dense range in $G$.  For every $\gamma\in\widehat G$, the composition $t\mapsto\gamma(\T_t|_{\Rr})$ is a continuous character of $\R$, and hence equals $e^{i\lambda t}$ for some $\lambda\in\R$.  Thus $P_\gamma(\Rr)\subseteq\A_\lambda$.

Conversely, if $a\in\A_\lambda$, then its orbit lies in the compact circle $\{za:|z|=1\}$.  Decompose $a=P(a)+(\id-P)(a)$.  The second term is again a $\lambda$-eigenoperator, so its norm is constant under the semigroup, but it lies in the stable kernel and therefore tends to zero.  Hence $(\id-P)(a)=0$ and $a\in\Rr$.
\end{proof}

\begin{corollary}\label{cor:scalar-relaxation}
Under the assumptions of Theorem~\ref{thm:jdlg}, the following are equivalent:
\begin{enumerate}[label=\upshape(\roman*)]
  \item $\Rr=\C1$;
  \item $\Fix(\T)=\C1$ and $\A_\lambda=\{0\}$ for every $\lambda\neq0$;
  \item there is a unique invariant state $\pi$ and
  \[
     \|\T_t(a)-\pi(a)1\|\longrightarrow0
     \qquad(a\in\A).
  \]
\end{enumerate}
\end{corollary}

\begin{proof}
The equivalence of \upshape(i) and \upshape(ii) follows from Theorem~\ref{thm:spectral-synthesis}.  If $\Rr=\C1$, write $P(a)=\pi(a)1$; positivity and unitality make $\pi$ a state, and stability of the kernel gives \upshape(iii).  Any invariant state equals its composition with $P$ by Theorem~\ref{thm:invariant-states}, so it equals $\pi$.  Finally, let $r\in\Rr$ and assume \upshape(iii).  The maps $\T_t|_{\Rr}$ are isometries for the Choi--Effros $C^*$-norm, which is the ambient norm.  Since $\pi$ is invariant,
\[
   \|r-\pi(r)1\|
   =\|\T_t(r-\pi(r)1)\|
   \longrightarrow0.
\]
Thus $r=\pi(r)1$ and $\Rr=\C1$.
\end{proof}

\section{The multiplicative reversible core}\label{sec:core}

For a unital completely positive map $P\colon\A\to\A$, its multiplicative domain is
\[
\begin{split}
   \MD(P):=\{a\in\A:
   &P(a^*a)=P(a)^*P(a),\\
   &P(aa^*)=P(a)P(a)^*\}.
\end{split}
\]
It is a $C^*$-subalgebra, and
\[
   P(ax)=P(a)P(x),
   \qquad
   P(xa)=P(x)P(a)
\]
whenever $a\in\MD(P)$ and $x\in\A$.  Positive projections and the operator-algebraic structure of their ranges have a substantial history; compare \cite{ChoiEffros1977,EffrosStormer1979}.  The next result isolates the part of a completely positive range on which the Choi--Effros and ambient products already agree.

\begin{theorem}\label{thm:multiplicative-core}
Let $P\colon\A\to\A$ be a unital completely positive projection, and put $\Rr=P(\A)$.  Then
\[
   \Rr_{\mathrm{amb}}:=\Rr\cap\MD(P)
   =\{a\in\Rr:a^*a,aa^*\in\Rr\}
\]
is a unital $C^*$-subalgebra of $\A$.  It is the largest ambient $C^*$-subalgebra contained in $\Rr$, and $P$ is $\Rr_{\mathrm{amb}}$-bimodular.  The following are equivalent:
\begin{enumerate}[label=\upshape(\roman*)]
  \item $\Rr_{\mathrm{amb}}=\Rr$;
  \item $\Rr$ is a $C^*$-subalgebra of $\A$ for the ambient product;
  \item $P$ is a conditional expectation onto $\Rr$.
\end{enumerate}
\end{theorem}

\begin{proof}
For $a\in\Rr$, one has $P(a)=a$.  Therefore
\[
   a\in\MD(P)
   \quad\Longleftrightarrow\quad
   P(a^*a)=a^*a\text{ and }P(aa^*)=aa^*,
\]
which is the displayed characterisation.  The intersection of the operator system $\Rr$ with the $C^*$-algebra $\MD(P)$ is linear and self-adjoint.  If $a,b\in\Rr_{\mathrm{amb}}$, then $ab\in\MD(P)$ and
\[
   P(ab)=P(a)P(b)=ab,
\]
so $ab\in\Rr_{\mathrm{amb}}$.  Thus it is a $C^*$-subalgebra.  Bimodularity follows from the defining property of the multiplicative domain.

Let $\B\subseteq\Rr$ be any ambient $C^*$-subalgebra.  If $b\in\B$, then $b^*b,bb^*\in\B\subseteq\Rr$, so the displayed characterisation gives $b\in\MD(P)$.  Hence $\B\subseteq\Rr_{\mathrm{amb}}$, proving maximality.

The equivalence of \upshape(i) and \upshape(ii) follows from maximality.  If \upshape(ii) holds, $P$ is a norm-one projection onto a $C^*$-subalgebra, and Tomiyama's theorem makes it a conditional expectation \cite{Tomiyama1957}.  The converse is immediate.
\end{proof}

\begin{corollary}\label{cor:faithful-family}
Let $P$ be the reversible projection of Theorem~\ref{thm:jdlg}.  The following are equivalent:
\begin{enumerate}[label=\upshape(\roman*)]
  \item the invariant states form a faithful family;
  \item $P$ is faithful, in the sense that $x\in\A_+$ and $P(x)=0$ imply $x=0$.
\end{enumerate}
If these conditions hold, then
\[
   \Rr_{\mathrm{amb}}=\Rr,
\]
so $\Rr$ is an ambient $C^*$-subalgebra and $P$ is a conditional expectation.
\end{corollary}

\begin{proof}
Assume \upshape(i) and let $x\in\A_+$ satisfy $P(x)=0$.  Every invariant state $\varphi$ satisfies $\varphi=\varphi\circ P$ by Theorem~\ref{thm:invariant-states}, so $\varphi(x)=0$.  Faithfulness of the family gives $x=0$.

Conversely, suppose that $P$ is faithful and let $0\neq x\in\A_+$.  Then $P(x)>0$ in the Choi--Effros algebra $\Rr$.  Choose a state $\psi_0\in\Sst(\Rr)$ with $\psi_0(P(x))>0$ and average it over the compact group:
\[
   \psi:=\int_G \psi_0\circ g\,dm_G(g).
\]
The integrand is non-negative on $P(x)$ and is strictly positive at the identity; continuity therefore gives $\psi(P(x))>0$.  The state $\varphi:=\psi\circ P$ is invariant by Theorem~\ref{thm:invariant-states}, and $\varphi(x)>0$.  Thus the invariant states form a faithful family.

Now fix $a\in\Rr$.  The Schwarz defect
\[
   d:=P(a^*a)-a^*a
\]
is positive.  Every invariant state annihilates $d$, because $\varphi\circ P=\varphi$.  Faithfulness of the family gives $d=0$.  Applying the same argument to $a^*$ yields $P(aa^*)=aa^*$, so $a\in\MD(P)$.  Theorem~\ref{thm:multiplicative-core} completes the proof.
\end{proof}

\begin{remark}
A single faithful invariant state is a special case of Corollary~\ref{cor:faithful-family}.  The proof is the familiar Schwarz-defect argument, but the faithful-family formulation is strictly more flexible and identifies the precise separation property needed.
\end{remark}

\begin{example}\label{ex:m3-core}
Let $\B$ be a unital $C^*$-algebra and identify $M_n(\B)$ with $M_n(\C)\otimes\B$.  Define
\[
   E_2\colon M_2(\B)\longrightarrow\B,
   \qquad
   E_2((b_{ij})):=\frac{b_{11}+b_{22}}2.
\]
Let $q\colon M_3(\B)\to M_2(\B)$ be compression to the upper-left corner and put
\[
   j(Y):=
   \begin{pmatrix}
      Y&0\\
      0&E_2(Y)
   \end{pmatrix}
   \qquad(Y\in M_2(\B)).
\]
Both maps are unital and completely positive, and $qj=\id$.  Hence
\[
   P:=jq\colon M_3(\B)\longrightarrow M_3(\B)
\]
is a unital completely positive projection.  The semigroup
\[
   \T_t:=e^{-t}\id+(1-e^{-t})P
\]
has compact orbits and reversible projection $P$.  Its Choi--Effros range is $*$-isomorphic to $M_2(\B)$, but it is not closed under the ambient product.  More precisely,
\[
   P(M_3(\B))\cap\MD(P)=1_3\otimes\B.
\]
No invariant state is faithful.  In particular, choosing $\B$ infinite dimensional gives infinite-dimensional examples in which the Choi--Effros and ambient products differ.
\end{example}

\begin{proof}
The semigroup law follows from $P^2=P$, every $\T_t$ is a convex combination of UCP maps, and $\T_t(x)\to P(x)$ in norm.  For $Y,Z\in M_2(\B)$,
\[
   j(Y)\circ j(Z)
   =P(j(Y)j(Z))
   =j(YZ),
\]
so $j$ identifies the Choi--Effros range with $M_2(\B)$.  The range is not an ambient algebra: for $Y=e_{11}\otimes1_{\B}$,
\[
   j(Y)^2=
   \begin{pmatrix}
      Y&0\\0&\tfrac14 1_{\B}
   \end{pmatrix}
   \neq j(Y^2).
\]

Let $x=j(Y)$.  A direct calculation gives
\[
   P(x^*x)=j(Y^*Y),
   \qquad
   x^*x=
   \begin{pmatrix}
      Y^*Y&0\\0&E_2(Y)^*E_2(Y)
   \end{pmatrix}.
\]
Thus $P(x^*x)=x^*x$ exactly when
\[
   E_2(Y^*Y)=E_2(Y)^*E_2(Y),
\]
and the analogous condition for $xx^*$ is obtained by replacing $Y$ with $Y^*$.  The map $E_2$ is faithful and $\B$-bimodular.  If the displayed equality holds and
\[
   Z:=Y-1_2\otimes E_2(Y),
\]
then $E_2(Z^*Z)=0$, hence $Z=0$.  Conversely, equality is immediate for $Y=1_2\otimes b$.  Therefore
\[
   P(M_3(\B))\cap\MD(P)
   =j(1_2\otimes\B)
   =1_3\otimes\B.
\]

Finally, an invariant state $\varphi$ satisfies $\varphi=\varphi\circ P$.  Since
\[
   P(e_{33}\otimes1_{\B})=0,
\]
every invariant state vanishes on the non-zero positive element $e_{33}\otimes1_{\B}$ and is therefore non-faithful.
\end{proof}

\section{Lip-norm criteria and metric consequences}\label{sec:lip}

\begin{definition}
Let $\A$ be a unital $C^*$-algebra and let $\Ldom\subseteq\A_{\mathrm{sa}}$ be a norm-dense unital real linear subspace. A seminorm $L\colon \Ldom\to[0,\infty)$ is called a \emph{Lip-norm} if:
\begin{enumerate}[label=\upshape(\roman*)]
    \item $L(a)=0$ if and only if $a\in\R\,1$,
    \item the Monge--Kantorovich metric
    \[
    d_L(\varphi,\psi):=\sup\bigl\{ |\varphi(a)-\psi(a)| : a\in\Ldom,\ L(a)\leqslant 1 \bigr\}
    \]
    induces the weak-$^*$ topology on $\Sst(\A)$.
\end{enumerate}
\end{definition}

\begin{remark}
This is the standard compact quantum metric-space framework of Rieffel \cite{Rieffel1999}. The only point we need is that weak-$^*$ neighbourhoods of a state may be described by small $d_L$-balls.
\end{remark}

\begin{lemma}\label{lem:mk-estimate}
Let $L$ be a Lip-norm on $\Ldom$. Then for every $a\in\Ldom$ and every $\varphi,\psi\in\Sst(\A)$ one has
\[
|\varphi(a)-\psi(a)|\leqslant L(a)\,d_L(\varphi,\psi).
\]
\end{lemma}

\begin{proof}
If $L(a)=0$, then $a$ is scalar and both sides vanish. If $L(a)>0$, apply the definition of $d_L$ to $a/L(a)$.
\end{proof}

\begin{proposition}\label{prop:uniform-lipschitz}
Let $\A\subseteq\M$ be a Feller algebra for $(\T_t)$, let $L$ be a Lip-norm on $\Ldom$, and assume that $\T_t(\Ldom)\subseteq\Ldom$ for every $t\geqslant 0$. If there exists $C\geqslant 0$ such that
\begin{equation}\label{eq:uniform-lip}
L(\T_t(a))\leqslant C\,L(a), \qquad a\in\Ldom,\ t\geqslant 0,
\end{equation}
then $(\T_t)$ has the $e$-property on $\Ldom$.
\end{proposition}

\begin{proof}
Fix $a\in\Ldom$ and $\varphi,\psi\in\Sst(\A)$. By Lemma~\ref{lem:mk-estimate} and \eqref{eq:uniform-lip},
\[
|\psi(\T_t(a)) - \varphi(\T_t(a))| \leqslant L(\T_t(a))\,d_L(\psi,\varphi) \leqslant C\,L(a)\,d_L(\psi,\varphi)
\]
for all $t\geqslant 0$. Since $d_L$ induces the weak-$^*$ topology, the required equicontinuity follows.
\end{proof}

\begin{corollary}\label{cor:eventual-lipschitz}
Let $\A\subseteq\M$ be a Feller algebra for $(\T_t)$, let $L$ be a Lip-norm on $\Ldom$, and assume that $\T_t(\Ldom)\subseteq\Ldom$ for every $t\geqslant 0$. If there exist $t_0\geqslant 0$ and $C\geqslant 0$ such that
\[
L(\T_t(a))\leqslant C\,L(a), \qquad a\in\Ldom,\ t\geqslant t_0,
\]
then $(\T_t)$ has the eventual $e$-property on $\Ldom$, and therefore the $e$-property on $\Ldom$.
\end{corollary}

\begin{proof}
The first claim is proved exactly as in Proposition~\ref{prop:uniform-lipschitz}, but only for $t\geqslant t_0$. The second claim follows from Theorem~\ref{thm:eventual-implies-e}.
\end{proof}

\begin{corollary}\label{cor:derivation}
Let $\A\subseteq\M$ be a Feller algebra for $(\T_t)$, and let $\delta\colon \dom(\delta)\to E$ be a closed $*$-derivation from a dense unital $*$-subalgebra $\dom(\delta)\subseteq\A$ into a Banach $\A$-bimodule $E$. Assume that $\T_t(\dom(\delta))\subseteq\dom(\delta)$ for every $t\geqslant 0$. Put
\[
L(a):=\|\delta(a)\|_E, \qquad a\in \dom(\delta)_{\mathrm{sa}},
\]
and assume that $L$ is a Lip-norm on $\dom(\delta)_{\mathrm{sa}}$. If
\begin{equation}\label{eq:derivation-bound}
\|\delta(\T_t(a))\|_E \leqslant C\,\|\delta(a)\|_E, \qquad a\in\dom(\delta),\ t\geqslant 0,
\end{equation}
then $(\T_t)$ has the $e$-property on $\dom(\delta)_{\mathrm{sa}}$.

More generally, if \eqref{eq:derivation-bound} is assumed only for $t\geqslant t_0$, then $(\T_t)$ has the eventual $e$-property and hence the $e$-property.
\end{corollary}

\begin{proof}
Apply Proposition~\ref{prop:uniform-lipschitz} and Corollary~\ref{cor:eventual-lipschitz} with $L(a)=\|\delta(a)\|_E$.
\end{proof}

\begin{remark}
Here ``closed'' means that the graph of $\delta$ is closed in $\A\times E$, or equivalently that $\dom(\delta)$ is complete for the graph norm
\[
\|x\|_{\delta}:=\|x\|+\|\delta(x)\|_E.
\]
\end{remark}

\begin{theorem}\label{thm:strict-contraction}
Let $\A\subseteq\M$ be a Feller algebra for $(\T_t)$, let $L$ be a Lip-norm on $\Ldom$, and suppose that $\T_t(\Ldom)\subseteq\Ldom$ for all $t\geqslant 0$. Assume that there exist $t_0>0$ and $\kappa\in[0,1)$ such that
\begin{equation}\label{eq:strict-contraction}
L(\T_{t_0}(a))\leqslant \kappa\,L(a), \qquad a\in\Ldom.
\end{equation}
Then the map $\Phi\colon \Sst(\A)\to\Sst(\A)$, defined by $\Phi(\varphi)=\varphi\circ\T_{t_0}$, is a strict contraction on the compact metric space $(\Sst(\A),d_L)$. Consequently:
\begin{enumerate}[label=\upshape(\roman*)]
    \item there exists a unique state $\varphi_*\in\Sst(\A)$ such that $\varphi_*\circ\T_{t_0}=\varphi_*$,
    \item $\varphi_*$ is invariant for the whole semigroup,
    \item for every $\varphi\in\Sst(\A)$ and every $n\in\N$,
    \[
    d_L\bigl(\varphi\circ\T_{n t_0},\varphi_*\bigr)\leqslant \kappa^n d_L(\varphi,\varphi_*).
    \]
\end{enumerate}
If, in addition, there exists $C\geqslant 0$ such that $L(\T_t(a))\leqslant C L(a)$ for all $a\in\Ldom$ and all $0\leqslant t\leqslant t_0$, then for every $t\geqslant 0$,
\[
 d_L\bigl(\varphi\circ\T_t,\varphi_*\bigr)\leqslant C\,\kappa^{\lfloor t/t_0\rfloor} d_L(\varphi,\varphi_*).
\]
\end{theorem}

\begin{proof}
For $\varphi,\psi\in\Sst(\A)$ and $a\in\Ldom$ with $L(a)\leqslant 1$, Lemma~\ref{lem:mk-estimate} and \eqref{eq:strict-contraction} give
\[
|\varphi(\T_{t_0}(a)) - \psi(\T_{t_0}(a))| \leqslant L(\T_{t_0}(a))\, d_L(\varphi,\psi) \leqslant \kappa\, d_L(\varphi,\psi).
\]
Taking the supremum over $a$ yields
\[
d_L\bigl(\Phi(\varphi),\Phi(\psi)\bigr)\leqslant \kappa\, d_L(\varphi,\psi).
\]
Since $(\Sst(\A),d_L)$ is compact and therefore complete, Banach's fixed-point theorem gives a unique fixed point $\varphi_*$ of $\Phi$. This proves \upshape(i) and \upshape(iii).

To prove \upshape(ii), fix $s\geqslant 0$. Since the dual maps commute,
\[
\bigl(\varphi_*\circ\T_s\bigr)\circ\T_{t_0} = \varphi_*\circ\T_{s+t_0} = \varphi_*\circ\T_{t_0+s} = \bigl(\varphi_*\circ\T_{t_0}\bigr)\circ\T_s = \varphi_*\circ\T_s.
\]
Thus $\varphi_*\circ\T_s$ is another fixed point of $\Phi$, hence equals $\varphi_*$ by uniqueness.

For the final estimate, write $t=n t_0+r$ with $n=\lfloor t/t_0\rfloor$ and $0\leqslant r<t_0$. Then
\begin{align*}
 d_L\bigl(\varphi\circ\T_t,\varphi_*\bigr)
 &= d_L\bigl((\varphi\circ\T_r)\circ\T_{n t_0},\varphi_*\bigr)
 \leqslant \kappa^n d_L\bigl(\varphi\circ\T_r,\varphi_*\bigr) \\
 &= \kappa^n d_L\bigl(\varphi\circ\T_r,\varphi_*\circ\T_r\bigr)
 \leqslant C\,\kappa^n d_L(\varphi,\varphi_*),
\end{align*}
as claimed.
\end{proof}

\begin{remark}
Theorem~\ref{thm:strict-contraction} is the metric-state-space version of the classical Wasserstein contraction principle. It is also one of the cleanest places where the commutative and the non-commutative theories genuinely coincide.
\end{remark}

\begin{corollary}\label{cor:strict-lip-norm}
Assume the hypotheses of Theorem~\ref{thm:strict-contraction}.  Suppose also that
\[
\diam(\Sst(\A),d_L)<\infty
\]
and that
\[
L(\T_t(a))\leqslant C L(a)
\qquad(a\in\Ldom,\ 0\leqslant t\leqslant t_0).
\]
Then, for the unique invariant state $\varphi_*$ and every $a\in\Ldom$,
\[
\|\T_t(a)-\varphi_*(a)1\|
\leqslant
C\,\diam(\Sst(\A),d_L)\,
\kappa^{\lfloor t/t_0\rfloor}L(a).
\]
\end{corollary}

\begin{proof}
Use the dual contraction estimate from Theorem~\ref{thm:strict-contraction}, test it on $a/L(a)$, and take the supremum over the state space.
\end{proof}

\section{Doeblin minorisation and strict-contraction maps}\label{sec:doeblin}

The global classical Doeblin condition says that one step of the dynamics dominates a fixed reference probability measure.  In the present setting it is a strong sufficient condition for the $e$-property and belongs beside non-commutative Dobrushin and quantum Doeblin coefficients \cite{GaubertQu2015,GeorgeHircheNuradhaWilde2026}. In operator-algebraic language the exact rank-one analogue is domination by the unital positive rank-one map $a\mapsto \sigma(a)1$. There is, however, a slightly more flexible formulation which turns out to be just as easy to prove and much better adapted to non-commutative applications.

\begin{theorem}\label{thm:general-doeblin}
Let $T\colon \A\to\A$ be a unital positive map on a unital $C^*$-algebra $\A$. Assume that there exist a unital positive map $\Psi\colon \A\to\A$, a constant $\lambda\in(0,1]$, and a number $\beta\in[0,1)$ such that
\begin{equation}\label{eq:minorisation-general}
T(a)\geqslant \lambda\,\Psi(a), \qquad a\in\A_+,
\end{equation}
and
\begin{equation}\label{eq:psi-contraction}
\|\varphi\circ\Psi - \psi\circ\Psi\| \leqslant \beta\,\|\varphi-\psi\|,
\qquad \varphi,\psi\in\Sst(\A).
\end{equation}
Put
\[
q:=1-\lambda(1-\beta)\in[0,1).
\]
Then:
\begin{enumerate}[label=\upshape(\roman*)]
    \item the dual map $T^*\colon \Sst(\A)\to\Sst(\A)$ is a strict contraction with Lipschitz constant at most $q$ in the norm of $\A^*$,
    \item there exists a unique invariant state $\varphi_*\in\Sst(\A)$ for $T$,
    \item for every $\varphi,\psi\in\Sst(\A)$ and every $n\in\N$,
    \[
    \bigl\|\varphi\circ T^n - \psi\circ T^n\bigr\| \leqslant q^n\,\|\varphi-\psi\|,
    \]
    \item for every $a\in\A_{\mathrm{sa}}$ and every $n\in\N$,
    \[
    \sup_{\varphi,\psi\in\Sst(\A)}\bigl|\varphi(T^n(a)) - \psi(T^n(a))\bigr| \leqslant 2\|a\|q^n.
    \]
\end{enumerate}
\end{theorem}

\begin{proof}
If $\lambda=1$, then $T(a)\geqslant \Psi(a)$ for all $a\geqslant 0$. Applying this to $1-a$ whenever $0\leqslant a\leqslant 1$ yields
\[
T(1-a)\geqslant \Psi(1-a),
\]
that is,
\[
-T(a)\geqslant -\Psi(a).
\]
Hence $T(a)\leqslant \Psi(a)$ on the positive unit ball. Combined with the original inequality, this gives $T=\Psi$ on the positive unit ball and therefore everywhere by linearity. In this case \upshape(i)--\upshape(iv) follow directly from \eqref{eq:psi-contraction}. Assume therefore that $0<\lambda<1$.

Define
\[
Q:=\frac{T-\lambda\Psi}{1-\lambda}.
\]
By \eqref{eq:minorisation-general}, the map $Q$ is positive on $\A_+$. Since $T(1)=\Psi(1)=1$, we also have $Q(1)=1$. Hence $Q$ is a positive unital map. Therefore its dual $Q^*$ preserves the state space and is contractive in norm. For states $\varphi,\psi$,
\[
\varphi\circ T - \psi\circ T = \lambda(\varphi\circ\Psi - \psi\circ\Psi) + (1-\lambda)(\varphi\circ Q - \psi\circ Q).
\]
Using \eqref{eq:psi-contraction} and the contractivity of $Q^*$, we obtain
\begin{align*}
\|\varphi\circ T - \psi\circ T\|
&\leqslant \lambda\,\|\varphi\circ\Psi - \psi\circ\Psi\| + (1-\lambda)\,\|\varphi\circ Q - \psi\circ Q\| \\
&\leqslant \lambda\beta\,\|\varphi-\psi\| + (1-\lambda)\,\|\varphi-\psi\| \\
&= q\,\|\varphi-\psi\|.
\end{align*}
This proves \upshape(i).

The state space is norm closed in $\A^*$ and therefore complete in the dual norm. Banach's fixed-point theorem now gives a unique fixed point $\varphi_*$ of $T^*$, proving \upshape(ii). Iterating the contraction estimate from \upshape(i) gives \upshape(iii). Finally, for $a\in\A_{\mathrm{sa}}$,
\begin{align*}
\bigl|\varphi(T^n(a)) - \psi(T^n(a))\bigr|
&\leqslant \|a\|\,\bigl\|\varphi\circ T^n - \psi\circ T^n\bigr\| \\
&\leqslant \|a\| q^n \|\varphi-\psi\| \\
&\leqslant 2\|a\|q^n.
\end{align*}
This proves \upshape(iv).
\end{proof}

\begin{corollary}\label{cor:rank-one-doeblin}
In Theorem~\ref{thm:general-doeblin}, if $\Psi(a)=\sigma(a)1$ for some state $\sigma\in\Sst(\A)$, then one may take $\beta=0$ and hence $q=1-\lambda$. This is the exact non-commutative analogue of the global classical Doeblin minorisation.
\end{corollary}

\begin{proof}
For the rank-one map $\Psi(a)=\sigma(a)1$ one has $\varphi\circ\Psi=\sigma=\psi\circ\Psi$ for all states $\varphi,\psi$, so \eqref{eq:psi-contraction} holds with $\beta=0$. Apply Theorem~\ref{thm:general-doeblin}.
\end{proof}

\begin{corollary}\label{cor:doeblin-semigroup}
Let $(\T_t)_{t\geqslant 0}$ be a quantum Markov semigroup on $\M$, and let $\A\subseteq\M$ be a Feller algebra. Assume that there exist $t_0>0$, a unital positive map $\Psi\colon \A\to\A$, and numbers $\lambda\in(0,1]$, $\beta\in[0,1)$ such that
\[
\T_{t_0}(a)\geqslant \lambda\,\Psi(a), \qquad a\in\A_+,
\]
and \eqref{eq:psi-contraction} holds for $\Psi$. Put $q=1-\lambda(1-\beta)$. Then:
\begin{enumerate}[label=\upshape(\roman*)]
    \item there exists a unique state $\varphi_*\in\Sst(\A)$ invariant for the whole semigroup,
    \item for every $\varphi,\psi\in\Sst(\A)$ and every $t\geqslant 0$,
    \[
    \|\varphi\circ\T_t - \psi\circ\T_t\|\leqslant q^{\lfloor t/t_0\rfloor}\,\|\varphi-\psi\|,
    \]
    \item $(\T_t)$ has the eventual $e$-property on $\A_{\mathrm{sa}}$, in fact with the neighbourhood $U=\Sst(\A)$, and therefore it has the $e$-property on $\A_{\mathrm{sa}}$.
\end{enumerate}
\end{corollary}

\begin{proof}
Apply Theorem~\ref{thm:general-doeblin} to the skeleton map $T:=\T_{t_0}|_{\A}$. This yields a unique state $\varphi_*$ with $\varphi_*\circ\T_{t_0}=\varphi_*$. As in the standard semigroup argument,
\[
\bigl(\varphi_*\circ\T_s\bigr)\circ\T_{t_0}=\varphi_*\circ\T_{s+t_0}=\varphi_*\circ\T_{t_0+s}=\bigl(\varphi_*\circ\T_{t_0}\bigr)\circ\T_s=\varphi_*\circ\T_s,
\]
so uniqueness of the fixed point of $\T_{t_0}^*$ implies $\varphi_*\circ\T_s=\varphi_*$ for every $s\geqslant 0$. This proves \upshape(i).

For \upshape(ii), write $t=n t_0+r$ with $n=\lfloor t/t_0\rfloor$ and $0\leqslant r<t_0$. Then
\begin{align*}
\|\varphi\circ\T_t - \psi\circ\T_t\|
&=\| (\varphi\circ\T_r)\circ\T_{n t_0} - (\psi\circ\T_r)\circ\T_{n t_0}\| \\
&\leqslant q^n\,\|\varphi\circ\T_r - \psi\circ\T_r\| \\
&\leqslant q^n\,\|\varphi-\psi\|,
\end{align*}
because the dual action of every unital positive map is norm contractive in the dual norm.

To prove \upshape(iii), fix $a\in\A_{\mathrm{sa}}$ and $\varepsilon>0$. Choose $N\in\N$ such that $2\|a\|q^N<\varepsilon$, and put $\tau:=N t_0$. Then for every $\varphi,\psi\in\Sst(\A)$ and every $t\geqslant \tau$,
\[
|\varphi(\T_t(a)) - \psi(\T_t(a))|\leqslant \|a\|\,\|\varphi\circ\T_t - \psi\circ\T_t\|\leqslant 2\|a\|q^{\lfloor t/t_0\rfloor}\leqslant 2\|a\|q^N<\varepsilon.
\]
Thus the eventual $e$-property holds globally, with the whole state space as neighbourhood. Theorem~\ref{thm:eventual-implies-e} now gives the $e$-property.
\end{proof}

\begin{corollary}\label{cor:factorised-minorisation}
Let $\eta\colon \A\to M_m(\C)$, $\theta\colon M_m(\C)\to\A$, and $R\colon M_m(\C)\to M_m(\C)$ be unital completely positive maps, and put
\[
\Psi:=\theta\circ R\circ \eta.
\]
Assume that there exists $\beta\in[0,1)$ such that
\[
\|\rho\circ R - \omega\circ R\| \leqslant \beta\,\|\rho-\omega\|,
\qquad \rho,\omega\in\Sst(M_m(\C)).
\]
Then $\Psi$ satisfies \eqref{eq:psi-contraction} with the same constant $\beta$. Consequently, every semigroup minorisation of the form
\[
\T_{t_0}(a)\geqslant \lambda\,\Psi(a), \qquad a\in\A_+,
\]
falls under Corollary~\ref{cor:doeblin-semigroup}.
\end{corollary}

\begin{proof}
Let $\varphi,\psi\in\Sst(\A)$. Since $\eta$ and $\theta$ are unital positive maps, their duals are norm contractive on state spaces. Moreover, $\varphi\circ\theta$ and $\psi\circ\theta$ are states on $M_m(\C)$, so the contraction hypothesis for $R$ applies to them. Therefore
\begin{align*}
\|\varphi\circ\Psi - \psi\circ\Psi\|
&=\| (\varphi\circ\theta)\circ R\circ\eta - (\psi\circ\theta)\circ R\circ\eta\| \\
&\leqslant \|(\varphi\circ\theta)\circ R - (\psi\circ\theta)\circ R\| \\
&\leqslant \beta\,\|\varphi\circ\theta - \psi\circ\theta\| \\
&\leqslant \beta\,\|\varphi-\psi\|.
\end{align*}
Hence \eqref{eq:psi-contraction} holds for $\Psi$.

\end{proof}

\begin{proposition}\label{prop:composable-minorisations}
Let $S,T\colon \A\to\A$ be unital positive maps. Assume that there exist positive maps $\Psi,\Phi\colon \A\to\A$, constants $\lambda,\kappa\in(0,1]$, and $\beta\in[0,1)$ such that
\begin{equation}\label{eq:minorisation-composable-1}
T(a)\geqslant \lambda\,\Psi(a),
\qquad a\in\A_+,
\end{equation}
\begin{equation}\label{eq:minorisation-composable-2}
S(\Psi(a))\geqslant \kappa\,\Phi(a),
\qquad a\in\A_+,
\end{equation}
and
\begin{equation}\label{eq:phi-contraction}
\|\varphi\circ\Phi-\psi\circ\Phi\|\leqslant \beta\,\|\varphi-\psi\|,
\qquad \varphi,\psi\in\Sst(\A).
\end{equation}
If, moreover, $\Phi$ is unital, then
\[
ST(a)\geqslant \lambda\kappa\,\Phi(a),
\qquad a\in\A_+,
\]
and therefore $ST$ falls under Theorem~\ref{thm:general-doeblin} with the same constant $\beta$.
\end{proposition}

\begin{proof}
If $a\in\A_+$, then by positivity of $S$ and \eqref{eq:minorisation-composable-1},
\[
ST(a)=S(T(a))\geqslant \lambda\,S(\Psi(a)).
\]
Using \eqref{eq:minorisation-composable-2}, we obtain
\[
ST(a)\geqslant \lambda\kappa\,\Phi(a).
\]
Since $\Phi$ is unital and satisfies \eqref{eq:phi-contraction}, Theorem~\ref{thm:general-doeblin} applies to the map $ST$.
\end{proof}

\begin{corollary}\label{cor:positive-profile-doeblin}
Let $(\T_t)_{t\geqslant 0}$ be a quantum Markov semigroup on $\M$, and let $\A\subseteq\M$ be a Feller algebra. Assume that there exist an element $h\in\A_+$, a state $\sigma\in\Sst(\A)$, times $t_0,s>0$, and constants $\lambda,\kappa\in(0,1]$ such that
\begin{equation}\label{eq:profile-rankone-step}
\T_{t_0}(a)\geqslant \lambda\,\sigma(a)\,h,
\qquad a\in\A_+,
\end{equation}
and
\begin{equation}\label{eq:uniform-return-profile}
\T_s(h)\geqslant \kappa\,1.
\end{equation}
Then
\[
\T_{s+t_0}(a)\geqslant \lambda\kappa\,\sigma(a)\,1,
\qquad a\in\A_+.
\]
Consequently $(\T_t)$ falls under Corollary~\ref{cor:doeblin-semigroup}, hence has a unique invariant state on $\A$ and the $e$-property on $\A_{\mathrm{sa}}$.
\end{corollary}

\begin{proof}
Apply Proposition~\ref{prop:composable-minorisations} with
\[
T:=\T_{t_0}|_{\A},
\qquad
S:=\T_s|_{\A},
\qquad
\Psi(a):=\sigma(a)\,h,
\qquad
\Phi(a):=\sigma(a)\,1.
\]
The map $\Phi$ is unital and satisfies \eqref{eq:phi-contraction} with $\beta=0$, while \eqref{eq:minorisation-composable-2} is exactly the consequence of \eqref{eq:uniform-return-profile}
\[
S(\Psi(a))=\sigma(a)\,\T_s(h)\geqslant \kappa\,\sigma(a)\,1=\kappa\,\Phi(a).
\]
Hence
\[
\T_{s+t_0}(a)=ST(a)\geqslant \lambda\kappa\,\sigma(a)\,1,
\qquad a\in\A_+,
\]
as claimed.
\end{proof}

\begin{corollary}\label{cor:commuting-corner-doeblin}
Let $(\T_t)_{t\geqslant 0}$ be a quantum Markov semigroup on $\M$, and let $\A\subseteq\M$ be a Feller algebra. Assume that there exist a projection $p\in\A$, a state $\sigma\in\Sst(\A)$, times $t_0,s>0$, and constants $\lambda,\kappa\in(0,1]$ such that
\[
p\,\T_{t_0}(a)\,p\geqslant \lambda\,\sigma(a)\,p,
\qquad a\in\A_+,
\]
\[
\T_s(p)\geqslant \kappa\,1,
\]
and, in addition,
\[
[p,\T_{t_0}(a)]=0
\qquad (a\in\A_+).
\]
Then the conclusion of Corollary~\ref{cor:positive-profile-doeblin} holds.

In particular, the conclusion holds whenever $p\in Z(\A)$.
\end{corollary}

\begin{proof}
If $x\in\A_+$ commutes with $p$, then
\[
x=pxp+(1-p)x(1-p)\geqslant pxp.
\]
Applying this to $x=\T_{t_0}(a)$ gives
\[
\T_{t_0}(a)\geqslant p\,\T_{t_0}(a)\,p\geqslant \lambda\,\sigma(a)\,p,
\qquad a\in\A_+.
\]
Now apply Corollary~\ref{cor:positive-profile-doeblin}.
\end{proof}

\begin{remark}
The commutation hypothesis in Corollary~\ref{cor:commuting-corner-doeblin} requires every element of $\T_{t_0}(\A_+)$ to commute with $p$. It ensures that the corner lower bound extends to a global order lower bound, as in the commutative setting.
\end{remark}

\begin{corollary}\label{cor:mean-return-profile}
Let $(\T_t)_{t\geqslant 0}$ be a quantum Markov semigroup on $\M$, and let $\A\subseteq\M$ be a Feller algebra. Assume that there exist an element $h\in\A_+$, a state $\sigma\in\Sst(\A)$, a time $t_0>0$, a number $S>0$, and constants $\lambda,\kappa\in(0,1]$ such that
\[
\T_{t_0}(a)\geqslant \lambda\,\sigma(a)\,h,
\qquad a\in\A_+,
\]
and
\begin{equation}\label{eq:mean-return-profile}
\frac1S\int_0^S \T_s(h)\,ds\geqslant \kappa\,1.
\end{equation}
Then the averaged map
\[
R:=M_S\circ \T_{t_0}\colon \A\to\A
\]
satisfies the global rank-one Doeblin minorisation
\[
R(a)\geqslant \lambda\kappa\,\sigma(a)\,1,
\qquad a\in\A_+.
\]
Consequently $R$ has a unique invariant state $\varphi_*\in\Sst(\A)$, this state is invariant for the whole semigroup $(\T_t)$, the dual of $R$ contracts the whole state space exponentially in norm, and for every $a\in\A_{\mathrm{sa}}$ and every $n\in\N$,
\[
\|R^n(a)-\varphi_*(a)1\|\leqslant 2\|a\|(1-\lambda\kappa)^n.
\]
\end{corollary}

\begin{proof}
If $a\in\A_+$, then
\[
R(a)=M_S(\T_{t_0}(a))\geqslant \lambda\,\sigma(a)\,M_S(h)\geqslant \lambda\kappa\,\sigma(a)\,1.
\]
Thus $R$ falls under Corollary~\ref{cor:rank-one-doeblin}. Let $\varphi_*$ denote its unique invariant state.

For every $u\geqslant 0$, the maps $\T_u$, $M_S$, and $\T_{t_0}$ commute. Hence
\[
(\varphi_*\circ \T_u)\circ R
=
(\varphi_*\circ R)\circ \T_u
=
\varphi_*\circ \T_u.
\]
Thus $\varphi_*\circ \T_u$ is also $R$-invariant, so uniqueness gives $\varphi_*\circ \T_u=\varphi_*$.

Finally, if $a\in\A_{\mathrm{sa}}$, then $R^n(a)-\varphi_*(a)1$ is self-adjoint and $\varphi_*$ is $R$-invariant, so
\[
\|R^n(a)-\varphi_*(a)1\|
=
\sup_{\psi\in\Sst(\A)}
\bigl|\psi(R^n(a))-\varphi_*(R^n(a))\bigr|
\leqslant
2\|a\|(1-\lambda\kappa)^n
\]
by Theorem~\ref{thm:general-doeblin}.
\end{proof}

\begin{remark}
A corner estimate of the form
\[
p\,\T_{t_0}(a)\,p\geqslant \lambda\,\sigma(a)\,p
\]
does \emph{not} by itself imply the global lower bound
\[
\T_{t_0}(a)\geqslant \lambda\,\sigma(a)\,p,
\]
because for non-central projections one cannot infer $x\geqslant pxp$ from positivity of $x$. This is precisely why the observable-by-observable projection-face theorem of the companion paper \cite{KaniaNiwinskiProjectionFaces} is needed as the genuine local non-commutative substitute for the classical small-set peeling argument. By contrast, the global Doeblin mechanism above only needs a positive return profile $h\in\A_+$ satisfying \eqref{eq:profile-rankone-step} and \eqref{eq:uniform-return-profile}; projections enter there only as a particularly important special case.
\end{remark}

\begin{remark}
Corollary~\ref{cor:mean-return-profile} is the natural Ces\`aro version of Corollary~\ref{cor:positive-profile-doeblin}: mean return to the positive profile $h$ already gives a genuine Doeblin estimate for the one-step averaged map $R=M_S\circ\T_{t_0}$. What it controls directly is this averaged operator. To pass from such a statement to norm mean-ergodic convergence of the original continuous-time semigroup still requires additional input, and that is exactly the point where the Ces\`aro machinery of Section~\ref{sec:cesaro} re-enters.
\end{remark}

\begin{remark}
Corollary~\ref{cor:factorised-minorisation} explains the possible rôle of nuclearity or injectivity. Those approximation properties furnish only \emph{approximate} finite-dimensional unital completely positive factorisations of the identity map, and therefore a plentiful supply of candidate finite-dimensional smoothing channels $\Psi$; see, for example, \cite[Ch.~2]{BrownOzawa2008}. What they do \emph{not} provide by themselves is either an exact factorisation of the identity through a finite-dimensional algebra or the order lower bound $\T_{t_0}\geqslant \lambda\Psi$. Producing such an exact lower bound is the genuinely non-commutative part of the Doeblin programme.
\end{remark}
\begin{proposition}\label{prop:finite-dimensional-doeblin}
Let $T\colon M_n(\C)\to M_n(\C)$ be a unital completely positive map. Assume that $T(a)$ is invertible for every non-zero positive element $a\in M_n(\C)$. Then there exists $\lambda>0$ such that
\[
T(a)\geqslant \lambda\,\operatorname{tr}_n(a)1,
\qquad a\in M_n(\C)_+,
\]
where $\operatorname{tr}_n$ is the normalised trace on $M_n(\C)$. Consequently $T$ has a unique invariant state and its dual acts as a strict contraction on the state space.
\end{proposition}

\begin{proof}
Let
\[
K:=\{\rho\in M_n(\C)_+ : \operatorname{tr}_n(\rho)=1\}.
\]
This is a compact subset of $M_n(\C)$. The function
\[
\rho\longmapsto \lambda_{\min}(T(\rho))
\]
is continuous on $K$, where $\lambda_{\min}$ denotes the smallest eigenvalue. By assumption it is strictly positive on $K$, hence it attains a strictly positive minimum $\lambda$. Therefore
\[
T(\rho)\geqslant \lambda 1,\qquad \rho\in K.
\]
If $a\in M_n(\C)_+$ is non-zero, write $a=\operatorname{tr}_n(a)\rho$ with $\rho\in K$. Then
\[
T(a)=\operatorname{tr}_n(a)T(\rho)\geqslant \lambda\,\operatorname{tr}_n(a)1.
\]
The last assertion follows from Corollary~\ref{cor:rank-one-doeblin}.
\end{proof}

\begin{corollary}\label{cor:factorised-rankone}
Let $\eta\colon \A\to M_m(\C)$ and $\theta\colon M_m(\C)\to\A$ be unital completely positive maps, let $R\colon M_m(\C)\to M_m(\C)$ be a unital completely positive map that is positivity improving in the sense of Proposition~\ref{prop:finite-dimensional-doeblin}, and put
\[
\Psi:=\theta\circ R\circ \eta.
\]
Then there exists $\lambda_R>0$ such that
\[
\Psi(a)\geqslant \lambda_R\,(\operatorname{tr}_m\circ\eta)(a)\,1,
\qquad a\in\A_+.
\]
Consequently, every minorisation of the form
\[
\T_{t_0}(a)\geqslant \mu\,\Psi(a),
\qquad a\in\A_+,
\]
automatically improves to the rank-one Doeblin bound
\[
\T_{t_0}(a)\geqslant \mu\lambda_R\,(\operatorname{tr}_m\circ\eta)(a)\,1,
\qquad a\in\A_+.
\]
In particular, Corollary~\ref{cor:doeblin-semigroup} applies.
\end{corollary}

\begin{proof}
By Proposition~\ref{prop:finite-dimensional-doeblin}, there exists $\lambda_R>0$ such that
\[
R(b)\geqslant \lambda_R\,\operatorname{tr}_m(b)\,1,
\qquad b\in M_m(\C)_+.
\]
Therefore, for $a\in\A_+$,
\[
\Psi(a)=\theta(R(\eta(a)))\geqslant \lambda_R\,\operatorname{tr}_m(\eta(a))\,\theta(1)=\lambda_R\,(\operatorname{tr}_m\circ\eta)(a)\,1.
\]
The final assertion is immediate.
\end{proof}

\begin{proposition}\label{prop:measure-prepare-rankone}
Let $\A$ be a unital $C^*$-algebra.  Choose states
$\sigma_1,\dots,\sigma_m\in\Sst(\A)$ and positive elements
$h_1,\dots,h_m\in\A_+$ satisfying
\[
\sum_{j=1}^m h_j=1.
\]
Define the unital positive measure-prepare map
\[
\Psi(a):=\sum_{j=1}^m \sigma_j(a)h_j,
\qquad a\in\A.
\]
Assume that there are constants $r_j\geqslant 0$, not all zero, such that
\[
h_j\geqslant r_j1
\qquad (1\leqslant j\leqslant m).
\]
Put $r:=\sum_{j=1}^m r_j$ and
\[
\sigma:=\frac1r\sum_{j=1}^m r_j\sigma_j.
\]
Then
\[
\Psi(a)\geqslant r\,\sigma(a)1,
\qquad a\in\A_+.
\]
Consequently, every lower bound
\[
\T_{t_0}(a)\geqslant \mu\,\Psi(a),
\qquad a\in\A_+,
\]
upgrades to the rank-one Doeblin bound
\[
\T_{t_0}(a)\geqslant \mu r\,\sigma(a)1,
\qquad a\in\A_+.
\]
\end{proposition}

\begin{proof}
For $a\in\A_+$ one has
\[
\Psi(a)=\sum_{j=1}^m \sigma_j(a)h_j
\geqslant
\sum_{j=1}^m r_j\sigma_j(a)1
=
r\,\sigma(a)1.
\]
The final assertion is immediate. Notice that $r\leqslant1$, since $1=\sum_jh_j\geqslant r1$.
\end{proof}

\begin{remark}
In finite-dimensional quantum information language, maps of the form in Proposition~\ref{prop:measure-prepare-rankone} are the basic entanglement-breaking, or measure-prepare, channels. The proposition says that an exact lower bound by such a channel becomes a genuine rank-one Doeblin bound as soon as the prepared positive elements have a uniform scalar part. This is stronger than what nuclearity or injectivity provide by themselves: those approximation properties give approximate finite-dimensional factorisations of the identity, whereas the proposition requires an exact minorising map with a quantitative order lower bound.
\end{remark}

\begin{remark}
Proposition~\ref{prop:finite-dimensional-doeblin} gives the simplest finite-dimensional model for the non-commutative Doeblin condition. Corollary~\ref{cor:factorised-rankone} shows that exact finite-dimensional factorisations of a minorising channel, with a positivity-improving middle channel, upgrade factorised lower bounds to genuine rank-one Doeblin bounds. This clarifies the approximation-property discussion: exact factorisations of minorising channels can solve the lower-bound problem, while approximate factorisations coming from nuclearity or injectivity alone only produce candidate channels and not the required order inequalities. For standard finite-dimensional quantum-channel terminology and background, see for instance \cite{Watrous2018,Wolf2012}.
\end{remark}

\section{Quantitative transfer and noiseless dynamics}\label{sec:transfer}

For a bounded norm-continuous trajectory $x\colon[0,\infty)\to\A_{\mathrm{sa}}$, let $\chi(x)$ be the Hausdorff measure of non-compactness of its range and let $\eDef(x)$ be the equicontinuity defect of the associated family $\{\widehat{x(t)}:t\geqslant0\}$ on $\Sst(\A)$.

\begin{theorem}\label{thm:trajectory-transfer}
Let $x,y\colon[0,\infty)\to\A_{\mathrm{sa}}$ be bounded norm-continuous trajectories and put
\[
   \eta:=\limsup_{t\to\infty}\|x(t)-y(t)\|.
\]
Then
\[
   |\chi(x)-\chi(y)|\leqslant\eta,
   \qquad
   |\eDef(x)-\eDef(y)|\leqslant2\eta.
\]
If $\eta=0$, the norm omega-limit sets of $x$ and $y$ coincide.  In particular, one trajectory converges in norm if and only if the other does, and the limits are equal.
\end{theorem}

\begin{proof}
Fix $\varepsilon>0$.  There is $T\geqslant0$ such that
\[
   \|x(t)-y(t)\|\leqslant\eta+\varepsilon
   \qquad(t\geqslant T).
\]
The two tail ranges are therefore at Hausdorff distance at most $\eta+\varepsilon$.  The Hausdorff measure of non-compactness is Lipschitz with respect to Hausdorff distance, so their measures differ by at most $\eta+\varepsilon$.  Deleting the initial compact trajectories does not change either measure.

For the defect, on every $U\subseteq\Sst(\A)$ and $t\geqslant T$,
\[
   \osc_U(\widehat{x(t)})
   \leqslant\osc_U(\widehat{y(t)})+2(\eta+\varepsilon).
\]
Taking the relevant suprema and infima, and then interchanging $x$ and $y$, gives the second inequality after $\varepsilon\downarrow0$.  Again, initial norm-compact trajectory segments have defect zero.

If $\eta=0$ and $t_i\to\infty$, then $x(t_i)\to z$ in norm if and only if $y(t_i)\to z$, proving equality of omega-limit sets.
\end{proof}

\begin{corollary}\label{cor:expected-transfer}
Let $\B\subseteq\A$ be an invariant unital $C^*$-subalgebra and let $E\colon\A\to\B$ be a conditional expectation.  For $a=a^*\in\A$, put
\[
   \eta(a):=\limsup_{t\to\infty}\|\T_t(a-E(a))\|.
\]
Then
\[
\begin{split}
   |\chiT_{\T}(a)-\chiT_{\T|_{\B}}(E(a))|&\leqslant\eta(a),\\
   |\eDef_{\T}(a)-\eDef_{\T|_{\B}}(E(a))|&\leqslant2\eta(a).
\end{split}
\]
If $\eta(a)=0$, the two trajectories have the same norm omega-limit set.  Thus the $e$-property, norm relative compactness and norm convergence transfer in both directions.
\end{corollary}

\begin{proof}
Apply Theorem~\ref{thm:trajectory-transfer} to
\[
   x(t)=\T_t(a),
   \qquad
   y(t)=\T_t(E(a)).
\]
Their difference is $\T_t(a-E(a))$.  The norm measure of the second trajectory is the same whether it is computed in $\B$ or in $\A$, since the inclusion is isometric.  For the equicontinuity defect, let
\[
   r\colon\Sst(\A)\to\Sst(\B),\qquad r(\rho)=\rho|_{\B},
\]
and
\[
   i\colon\Sst(\B)\to\Sst(\A),\qquad i(\theta)=\theta\circ E.
\]
Then $r$ is a continuous surjection and $r\circ i=\id$.  The expectation family of $y$ on $\Sst(\A)$ is the pullback by $r$ of the corresponding family on $\Sst(\B)$, so pulling back finite covers gives one inequality between the defects.  Restricting to the compact section $i(\Sst(\B))$ gives the reverse inequality.  Hence the two defects agree, and the displayed estimates follow.
\end{proof}

\begin{example}\label{ex:irrational-dephasing}
Let $\omega_1,\omega_2\in\R$ with $\omega_1/\omega_2\notin\mathbb Q$, and define a flow on $\mathbb T^2$ by
\[
   \theta_t(z_1,z_2)
   =(e^{i\omega_1t}z_1,e^{i\omega_2t}z_2).
\]
Let $m$ be Haar measure on $\mathbb T^2$ and put
\[
   \M=L^\infty(\mathbb T^2,\mathfrak{B}(\mathbb{T}^2),m)\,\bar\otimes\,M_2(\C),
   \qquad
   \A=C(\mathbb T^2,M_2(\C)).
\]
Let $\D_2\subseteq M_2(\C)$ be the diagonal algebra and let $E_{\D_2}$ be the diagonal conditional expectation. Put
\[
   D_t=E_{\D_2}+e^{-t}(\id-E_{\D_2}).
\]
Define
\[
   (\T_tf)(z):=D_t(f(\theta_tz)).
\]
Then the formula extends to a weak-$^*$ continuous normal quantum Markov semigroup on $\M$, the algebra $\A$ is Feller, and the $e$-property holds on all of $\A_{\mathrm{sa}}$.  Its reversible projection is the fibrewise diagonal expectation
\[
   E\colon C(\mathbb T^2,M_2(\C))\longrightarrow C(\mathbb T^2,\D_2),
\]
and its stable sector is the off-diagonal kernel of $E$.  The reversible dynamics has compact but generally non-periodic orbits.  No time slice is a strict contraction on the whole state space in the dual norm.
\end{example}

\begin{proof}
The base maps form a strongly continuous group of $*$-automorphisms on $\A$ and a weak-$^*$ continuous normal group on $L^\infty(\mathbb T^2,m)$.  The family $(D_t)$ is a norm-continuous UCP semigroup because $E_{\D_2}$ is a UCP projection, and the two actions commute.

The diagonal algebra $\B=C(\mathbb T^2,\D_2)$ is invariant.  The closure of the one-parameter subgroup $(e^{i\omega_1t},e^{i\omega_2t})$ is the compact torus, so every orbit in $\B$ is relatively compact.  Moreover,
\[
   \|\T_t(f-Ef)\|=e^{-t}\|f-Ef\|.
\]
Corollary~\ref{cor:expected-transfer} yields the $e$-property on the whole algebra and identifies the reversible and stable sectors.

For example,
\[
   b(z_1,z_2)=
   \begin{pmatrix}\operatorname{Re}z_1&0\\0&\operatorname{Re}z_2\end{pmatrix}
\]
has a non-periodic compact orbit.  Indeed, a period $\tau>0$ would force both $e^{i\omega_1\tau}=1$ and $e^{i\omega_2\tau}=1$, contradicting irrationality of the ratio.

For $z\in\mathbb T^2$, let $\rho_z(f)=\operatorname{tr}_2(f(z))$.  Since $D_t$ preserves the trace,
\[
   \rho_z\circ\T_t=\rho_{\theta_tz}.
\]
Distinct points remain distinct under the flow, and a scalar function separating them shows that the dual-norm distance of the corresponding states remains two.  Hence no time slice has a strict dual-norm contraction coefficient.
\end{proof}

\section{Expected subalgebras and decoherence-free decompositions}\label{sec:expected-subalgebras}

We now turn to another source of equicontinuity: the existence of an invariant expected subalgebra. The guiding situation is that a conditional expectation $E\colon\A\to\B$ captures the diagonal or noiseless part of the dynamics, while the complementary sector $\ker E$ decays in norm. The $e$-property can then be lifted from $\B$ to the whole algebra.

\begin{definition}\label{def:cartan}
Let $\A$ be a unital $C^*$-algebra. A unital $C^*$-subalgebra $\D\subseteq\A$ is called a \emph{Cartan masa} if:
\begin{enumerate}[label=\upshape(\roman*)]
    \item $\D$ is maximal abelian in $\A$,
    \item the normaliser
    \[
    \mathcal N_\A(\D):=\{n\in\A: n\D n^*\subseteq \D\ \text{and}\ n^*\D n\subseteq \D\}
    \]
    generates $\A$ as a $C^*$-algebra,
    \item there exists a faithful conditional expectation $E_\D\colon \A\to\D$.
\end{enumerate}
This is the $C^*$-algebraic Cartan terminology used by Renault, building on Kumjian's $C^*$-diagonals \cite{Kumjian1986,Renault2008}.
\end{definition}

\begin{theorem}\label{thm:expected-subalgebra-lifting}
Let $(\T_t)_{t\geqslant 0}$ be a quantum Markov semigroup on $\M$, let $\A\subseteq\M$ be a Feller algebra, and let $\B\subseteq\A$ be a unital $C^*$-subalgebra admitting a conditional expectation $E\colon \A\to\B$. Assume that $\T_t(\B)\subseteq\B$ for every $t\geqslant 0$. Let $\Ldom\subseteq \A_{\mathrm{sa}}$ and $\Ldom_\B\subseteq \B_{\mathrm{sa}}$ be unital real linear test spaces such that
\[
E(\Ldom)\subseteq \Ldom_\B.
\]
Assume further that:
\begin{enumerate}[label=\upshape(\roman*)]
    \item the restricted semigroup $(\T_t|_\B)$ has the eventual $e$-property on $\Ldom_\B$,
    \item for every $a\in \Ldom$,
    \[
    \|\T_t(a-E(a))\|\longrightarrow 0
    \qquad (t\to\infty).
    \]
\end{enumerate}
Then $(\T_t)$ has the eventual $e$-property on $\Ldom$, and therefore the $e$-property on $\Ldom$.
\end{theorem}

\begin{proof}
Fix $a\in\Ldom$, $\varphi\in\Sst(\A)$, and $\varepsilon>0$. Put
\[
d:=E(a)\in \Ldom_\B,
\qquad
b:=a-d.
\]
Let
\[
r\colon \Sst(\A)\to \Sst(\B),
\qquad
r(\omega):=\omega|_\B,
\]
be the restriction map. Since $\B\subseteq \A$ is unital, $r$ is weak-$^*$ continuous.

Because $(\T_t|_\B)$ has the eventual $e$-property on $\Ldom_\B$, there exist a weak-$^*$ neighbourhood $V$ of $r(\varphi)$ in $\Sst(\B)$ and a time $T_0\geqslant 0$ such that
\[
\sup_{\rho\in V}\sup_{t\geqslant T_0}
\bigl|\rho(\T_t(d))-r(\varphi)(\T_t(d))\bigr|<\varepsilon/2.
\]
Put
\[
U:=r^{-1}(V),
\]
which is a weak-$^*$ neighbourhood of $\varphi$ in $\Sst(\A)$.

Choose $T_1\geqslant 0$ such that
\[
\|\T_t(b)\|<\varepsilon/4,
\qquad t\geqslant T_1.
\]
Then for every $\psi\in U$ and every $t\geqslant \max\{T_0,T_1\}$,
\begin{align*}
\bigl|\psi(\T_t(a))-\varphi(\T_t(a))\bigr|
&\leqslant \bigl|\psi(\T_t(d))-\varphi(\T_t(d))\bigr|
+\bigl|\psi(\T_t(b))\bigr|+\bigl|\varphi(\T_t(b))\bigr| \\
&= \bigl|r(\psi)(\T_t(d))-r(\varphi)(\T_t(d))\bigr|
+\bigl|\psi(\T_t(b))\bigr|+\bigl|\varphi(\T_t(b))\bigr| \\
&< \varepsilon/2+2\|\T_t(b)\| \\
&<\varepsilon.
\end{align*}
Thus $(\T_t)$ has the eventual $e$-property at $\varphi$ for the observable $a$. Since $a\in\Ldom$ and $\varphi\in\Sst(\A)$ were arbitrary, the eventual $e$-property holds on $\Ldom$. Theorem~\ref{thm:eventual-implies-e} now yields the $e$-property on $\Ldom$.
\end{proof}

\begin{remark}
The inclusion $E(\Ldom)\subseteq\Ldom_\B$ in Theorem~\ref{thm:expected-subalgebra-lifting} is the natural hypothesis for the proof. In applications where $\Ldom=\A_{\mathrm{sa}}$ and $\Ldom_\B=\B_{\mathrm{sa}}$, it is of course automatic and becomes equality after applying $E$ to $\B_{\mathrm{sa}}$.
\end{remark}

\begin{remark}
The proof of Theorem~\ref{thm:expected-subalgebra-lifting} uses only an invariant expected subalgebra. The advantage of Cartan masas and related structures is therefore structural rather than formal: they provide concrete settings in which one may hope to verify the off-diagonal decay
\[
\|\T_t(a-E(a))\|\to 0
\]
on a large test space. In this sense they furnish a non-commutative replacement for the classical idea of separating diagonal behaviour from off-diagonal transport.
\end{remark}

\begin{remark}\label{rem:cstar-level-semigroups}
Several examples below are written directly on a unital $C^*$-algebra.  Their normal second-dual extensions act on $\A^{**}$, and the original algebra is a Feller algebra.  They therefore remain within the von Neumann--Feller framework of Section~\ref{sec:framework}.
\end{remark}

\begin{corollary}\label{cor:expected-subalgebra-decoherence}
Let $\A$ be a unital $C^*$-algebra, let $\B\subseteq\A$ be a unital $C^*$-subalgebra, and let $E\colon \A\to\B$ be a conditional expectation. Then
\[
\T_t:=e^{-t}\id+(1-e^{-t})E,
\qquad t\geqslant 0,
\]
is a norm-continuous semigroup of unital completely positive maps on $\A$. Moreover, $\T_t$ has the $e$-property on $\A_{\mathrm{sa}}$.

In particular, Cartan masas provide one important special case.
\end{corollary}

\begin{proof}
Since $E$ is a unital completely positive idempotent, the semigroup law
\[
\T_s\T_t=e^{-(s+t)}\id+(1-e^{-(s+t)})E=\T_{s+t}
\]
is immediate. Norm continuity in $t$ is obvious.

Let $\Ldom_\B:=\B_{\mathrm{sa}}$ and $\Ldom:=\A_{\mathrm{sa}}$. The restriction of $(\T_t)$ to $\B$ is the identity semigroup, which has the $e$-property because for each $d\in \B_{\mathrm{sa}}$ the function $\Sst(\B)\ni \rho\mapsto \rho(d)$ is weak-$^*$ continuous. Also,
\[
\T_t(a-E(a))=e^{-t}(a-E(a))
\qquad (a\in \A_{\mathrm{sa}}),
\]
so the complementary part decays exponentially in norm. Theorem~\ref{thm:expected-subalgebra-lifting} therefore gives the $e$-property on $\A_{\mathrm{sa}}$.

\end{proof}

\begin{proposition}\label{prop:finite-dimensional-spectral-gap}
Let $\A$ be a finite-dimensional unital $C^*$-algebra, let $(\T_t)_{t\geqslant0}$ be a norm-continuous semigroup of unital completely positive maps on $\A$, and let $\omega$ be a faithful invariant state. Let $\B\subseteq\A$ be a unital $C^*$-subalgebra and let $E\colon\A\to\B$ be an $\omega$-preserving conditional expectation. Assume that
\begin{enumerate}[label=\upshape(\roman*)]
    \item $E\T_t=\T_tE=E$ for every $t\geqslant0$;
    \item there exists $\gamma>0$ such that, for the GNS norm $\|x\|_{2,\omega}:=\omega(x^*x)^{1/2}$,
    \[
    \|\T_t(b)\|_{2,\omega}\leqslant e^{-\gamma t}\|b\|_{2,\omega},
    \qquad b\in\ker E,
    \quad t\geqslant0.
    \]
\end{enumerate}
Then $(\T_t)$ has the $e$-property on $\A_{\mathrm{sa}}$.
\end{proposition}

\begin{proof}
Hypothesis \upshape(i) already forces the restricted dynamics on $\B$ to be trivial. Indeed, if $d\in\B$, then $d=E(d)$ and hence
\[
\T_t(d)=\T_tE(d)=E(d)=d,
\qquad t\geqslant0.
\]
Thus $\T_t|_\B=\id_\B$, so the restricted semigroup has the $e$-property on $\B_{\mathrm{sa}}$ automatically.

Since $\A$ is finite-dimensional and $\omega$ is faithful, the $C^*$-norm and the GNS norm are equivalent. Choose $C_\omega>0$ such that
\[
\|x\|\leqslant C_\omega\|x\|_{2,\omega},
\qquad x\in\A.
\]
If $a\in\A_{\mathrm{sa}}$, then $a-E(a)\in\ker E$, and the spectral-gap assumption gives
\[
\|\T_t(a-E(a))\|
\leqslant C_\omega e^{-\gamma t}\|a-E(a)\|_{2,\omega}
\longrightarrow0.
\]
The conclusion follows from Theorem~\ref{thm:expected-subalgebra-lifting}, with $\Ldom=\A_{\mathrm{sa}}$ and $\Ldom_\B=\B_{\mathrm{sa}}$.
\end{proof}

\begin{remark}
Detailed balance or KMS-symmetry is a standard way of producing the hypotheses of Proposition~\ref{prop:finite-dimensional-spectral-gap}: in finite dimensions the generator becomes self-adjoint on $L^2(\A,\omega)$, and a spectral gap on the orthogonal complement of the fixed sector gives the displayed $L^2$-decay. The proposition isolates the only additional ingredient needed for the present paper, namely the passage from $L^2$-decay to norm decay, which is automatic in finite dimensions and in infinite-dimensional settings under suitable hypercontractive or ultracontractive estimates. The broader detailed-balance and functional-inequality background includes the non-commutative Poincar\'e, logarithmic Sobolev, gradient-flow, and Gibbs-sampler literature \cite{KastoryanoTemme2013,KastoryanoBrandao2016,CarlenMaas2017,BrandaoKastoryano2019}.
\end{remark}

\begin{theorem}\label{thm:dfs-lifting}
Let $(\T_t)_{t\geqslant 0}$ be a quantum Markov semigroup on a von Neumann algebra $\M$, and let $\A\subseteq\M$ be a Feller algebra. Let $\Nc\subseteq \A$ be a unital $C^*$-subalgebra admitting a conditional expectation $E_{\Nc}\colon \A\to\Nc$. Assume that:
\begin{enumerate}[label=\upshape(\roman*)]
    \item there exists a strongly continuous one-parameter group $(\alpha_t)_{t\in\R}$ of $*$-automorphisms of $\Nc$ such that
    \[
    \T_t|_\Nc=\alpha_t
    \qquad
    (t\geqslant 0),
    \]
    \item there exists a unital real linear test space $\Ldom\subseteq \A_{\mathrm{sa}}$ such that
    \[
    E_{\Nc}(\Ldom)\subseteq \Ldom
    \]
    and
    \[
    \|\T_t(a-E_{\Nc}(a))\|\longrightarrow 0
    \qquad
    (t\to\infty)
    \]
    for every $a\in\Ldom$.
\end{enumerate}
Then the following are equivalent:
\begin{enumerate}[label=\upshape(\alph*)]
    \item $(\T_t)$ has the eventual $e$-property on $\Ldom$,
    \item $(\T_t)$ has the $e$-property on $\Ldom$,
    \item the semigroup $(\alpha_t)_{t\geqslant 0}$ has the $e$-property on $E_{\Nc}(\Ldom)\subseteq \Nc_{\mathrm{sa}}$.
\end{enumerate}
\end{theorem}

\begin{proof}
The implication \upshape(a)$\Rightarrow$\upshape(b) is exactly Theorem~\ref{thm:eventual-implies-e}. We prove \upshape(c)$\Rightarrow$\upshape(a). Fix $a\in\Ldom$, $\varphi\in\Sst(\A)$, and $\varepsilon>0$. Put
\[
d:=E_{\Nc}(a)\in E_{\Nc}(\Ldom)\subseteq \Ldom,
\qquad
b:=a-d.
\]
Because the semigroup $(\alpha_t)_{t\geqslant 0}$ has the $e$-property on $E_{\Nc}(\Ldom)$, there exists a weak-$^*$ neighbourhood $V$ of $\varphi|_\Nc$ in $\Sst(\Nc)$ such that
\[
\sup_{\rho\in V}\sup_{t\geqslant 0}
|\rho(\alpha_t(d))-(\varphi|_\Nc)(\alpha_t(d))|<\varepsilon/2.
\]
Let
\[
r\colon \Sst(\A)\to \Sst(\Nc),
\qquad
r(\omega)=\omega|_\Nc,
\]
which is weak-$^*$ continuous. Then $U:=r^{-1}(V)$ is a weak-$^*$ neighbourhood of $\varphi$ in $\Sst(\A)$. By hypothesis \upshape(ii), choose $T_0\geqslant 0$ such that
\[
\|\T_t(b)\|<\varepsilon/4,
\qquad t\geqslant T_0.
\]
For $\psi\in U$ and $t\geqslant T_0$ we obtain
\begin{align*}
|\psi(\T_t(a))-\varphi(\T_t(a))|
&\leqslant
|\psi(\T_t(d))-\varphi(\T_t(d))|
+|\psi(\T_t(b))|
+|\varphi(\T_t(b))| \\
&=
|r(\psi)(\alpha_t(d))-(\varphi|_\Nc)(\alpha_t(d))|
+|\psi(\T_t(b))|
+|\varphi(\T_t(b))| \\
&< \varepsilon/2+2\|\T_t(b)\| \\
&< \varepsilon.
\end{align*}
Thus $(\T_t)$ has the eventual $e$-property on $\Ldom$.

It remains to prove \upshape(b)$\Rightarrow$\upshape(c). By hypothesis \upshape(b), $(\T_t)$ has the $e$-property on $\Ldom$. Fix $d\in E_{\Nc}(\Ldom)$, $\rho\in\Sst(\Nc)$, and $\varepsilon>0$. Define
\[
i\colon \Sst(\Nc)\to \Sst(\A),
\qquad
i(\rho):=\rho\circ E_{\Nc}.
\]
Since $E_{\Nc}$ is unital and positive, $i$ is weak-$^*$ continuous. The $e$-property of $(\T_t)$ at the lifted state $i(\rho)$, evaluated on the observable $d\in \Ldom$, yields a weak-$^*$ neighbourhood $U$ of $i(\rho)$ in $\Sst(\A)$ such that
\[
\sup_{\psi\in U}\sup_{t\geqslant 0}
|\psi(\T_t(d))-i(\rho)(\T_t(d))|<\varepsilon.
\]
Put
\[
V:=i^{-1}(U),
\]
a weak-$^*$ neighbourhood of $\rho$ in $\Sst(\Nc)$. If $\omega\in V$, then $i(\omega)\in U$, and since $d\in \Nc$ and $\T_t|_\Nc=\alpha_t$ for $t\geqslant 0$,
\begin{align*}
|\omega(\alpha_t(d))-\rho(\alpha_t(d))|
&=
|i(\omega)(\T_t(d))-i(\rho)(\T_t(d))| \\
&< \varepsilon
\qquad (t\geqslant 0).
\end{align*}
Hence the semigroup $(\alpha_t)_{t\geqslant 0}$ has the $e$-property on $E_{\Nc}(\Ldom)$. This proves \upshape(c) and completes the equivalence.
\end{proof}

\begin{remark}
For a strongly continuous automorphism group $(\alpha_t)_{t\in\R}$ on a unital $C^*$-algebra $\Nc$ and a real linear subspace $\Ldom_0\subseteq \Nc_{\mathrm{sa}}$, the semigroup $(\alpha_t)_{t\geqslant 0}$ has the $e$-property on $\Ldom_0$ if and only if each positive orbit
\[
\{\alpha_t(d):t\geqslant 0\}
\qquad
(d\in\Ldom_0)
\]
is norm relatively compact. Indeed, norm relative compactness of the positive orbit gives the required weak-$^*$ equicontinuity by Lemma~\ref{lem:compact-uniform}. Conversely, if $(\alpha_t)_{t\geqslant 0}$ has the $e$-property on $\Ldom_0$, then the family of functions
\[
\widehat{\alpha_t(d)}\colon \Sst(\Nc)\to\R,
\qquad
\widehat{\alpha_t(d)}(\rho):=\rho(\alpha_t(d)),
\]
indexed by $t\geqslant 0$ is uniformly bounded and equicontinuous on the weak-$^*$ compact state space. Arzel\`a--Ascoli therefore makes it relatively compact in $C(\Sst(\Nc),\R)$, and since
\[
\|\widehat{x}-\widehat{y}\|_{\infty}=\|x-y\|
\qquad
(x,y\in \Nc_{\mathrm{sa}}),
\]
the positive orbit $\{\alpha_t(d):t\geqslant 0\}$ is norm relatively compact.

If one strengthens this to norm relative compactness of the full two-sided orbit $\{\alpha_t(d):t\in\R\}$, then the classical Bochner criterion identifies the orbit map $t\mapsto \alpha_t(d)$ with an almost periodic $\Nc$-valued function. In that stronger two-sided sense, Theorem~\ref{thm:dfs-lifting} recovers the familiar slogan that decoherence towards a decoherence-free subalgebra reduces the $e$-property to almost periodicity of the surviving noiseless dynamics.
\end{remark}

\begin{remark}\label{rem:finite-dimensional-trivial}
If $\A$ is finite-dimensional, then the weak-$^*$ and norm topologies on $\Sst(\A)$ coincide. By Remark~\ref{rem:norm-vs-weak}, every quantum Markov semigroup on $\A$ is therefore automatically norm-nonexpansive on states and hence has the $e$-property on every test space. The finite-dimensional examples below should thus be read as illustrations of the Cartan and decoherence-free mechanisms, not as genuinely nontrivial sources of the $e$-property itself.
\end{remark}

\begin{example}\label{ex:dfs-block-dephasing}
Let $\A=M_n(\C)$ and let $p_1,\dots,p_r$ be pairwise orthogonal projections with $\sum_{j=1}^r p_j=1$. Put
\[
\Nc:=\bigoplus_{j=1}^r p_jM_n(\C)p_j,
\qquad
E_{\Nc}(a):=\sum_{j=1}^r p_jap_j.
\]
Assume that at least one block $p_jM_n(\C)p_j$ has dimension greater than one, so that $\Nc$ is genuinely non-commutative. Choose $H=H^*\in\Nc$, let
\[
\alpha_t(a):=e^{itH}ae^{-itH},
\qquad a\in\Nc,
\]
and define
\[
\T_t(a):=
e^{itH}E_{\Nc}(a)e^{-itH}
+
e^{-t}e^{itH}(a-E_{\Nc}(a))e^{-itH},
\qquad
a\in M_n(\C),\ t\geqslant 0.
\]
Then $(\T_t)_{t\geqslant 0}$ is a norm-continuous quantum Markov semigroup on $M_n(\C)$, $E_{\Nc}\circ \T_t=\alpha_t\circ E_{\Nc}$, the restriction of $\T_t$ to $\Nc$ is the automorphism group $(\alpha_t)$, and
\[
\|\T_t(a-E_{\Nc}(a))\|=e^{-t}\|a-E_{\Nc}(a)\|
\qquad
(a\in M_n(\C)).
\]
Since $\Nc$ is finite-dimensional, every positive orbit (indeed every full orbit) of $(\alpha_t)$ is norm relatively compact. Hence the semigroup $(\alpha_t)_{t\geqslant 0}$ has the $e$-property on $\Nc_{\mathrm{sa}}$, and Theorem~\ref{thm:dfs-lifting} yields the $e$-property of $(\T_t)$ on $M_n(\C)_{\mathrm{sa}}$.
\end{example}

\begin{proof}
The semigroup law follows from the decomposition
\[
\T_t=\operatorname{Ad}(e^{itH})\circ\bigl(E_{\Nc}+e^{-t}(\id-E_{\Nc})\bigr)
\]
and the fact that $\operatorname{Ad}(e^{itH})$ commutes with $E_{\Nc}$ because $H\in\Nc$ and $E_{\Nc}$ is $\Nc$-bimodular. Complete positivity and unitality are immediate. The remaining assertions are exactly the formulas displayed above, and finite dimensionality gives norm relative compactness of the automorphism orbits.
\end{proof}

\begin{theorem}\label{thm:contractive-expected-subalgebra}
Let $\A$ be a unital $C^*$-algebra, let $\B\subseteq\A$ be a unital $C^*$-subalgebra, let $E\colon \A\to\B$ be a conditional expectation, and let $P\colon \A\to\A$ be a unital completely positive map such that
\begin{enumerate}[label=\upshape(\roman*)]
    \item $P(\B)\subseteq \B$ and $P|_\B=\id_\B$,
    \item $PE=EP=E$,
    \item there exists $\theta\in[0,1)$ such that
    \[
    \|P(b)\|\leqslant \theta\,\|b\|,
    \qquad b\in \ker E.
    \]
\end{enumerate}
Let $(\T_t)_{t\geqslant 0}$ be the uniformly continuous quantum Markov semigroup generated by $P-\id$, that is,
\[
\T_t=e^{t(P-\id)}=e^{-t}\sum_{n=0}^{\infty}\frac{t^n}{n!}P^n.
\]
Then:
\begin{enumerate}[label=\upshape(\roman*)]
    \item $\T_t(\B)=\B$ and $E\circ \T_t=\T_t\circ E=E$ for every $t\geqslant 0$,
    \item
    \[
    \|\T_t(a-E(a))\|\leqslant e^{-t(1-\theta)}\,\|a-E(a)\|,
    \qquad a\in\A,
    \]
    \item $(\T_t)$ has the $e$-property on $\A_{\mathrm{sa}}$.
\end{enumerate}
\end{theorem}

\begin{proof}
Because $P$ is unital and completely positive, each $P^n$ is unital and completely positive, and therefore so is $\T_t$ for every $t\geqslant 0$. Norm continuity is immediate from the series defining $\T_t$.

If $d\in\B$, then $P(d)=d$ by \upshape(i), whence $\T_t(d)=d$ for every $t\geqslant 0$. Thus $\T_t(\B)=\B$. Since $E$ commutes with $P$, it commutes with every polynomial in $P$, hence with every $\T_t$, and because $E\T_t|_\B=\id_\B$ one also has $E\circ\T_t=E$.

Now fix $b\in\ker E$. By \upshape(ii), $P(\ker E)\subseteq\ker E$, and by \upshape(iii),
\[
\|P^n(b)\|\leqslant \theta^n\|b\|,
\qquad n\in\N.
\]
Therefore
\[
\|\T_t(b)\|
\leqslant
e^{-t}\sum_{n=0}^{\infty}\frac{t^n}{n!}\|P^n(b)\|
\leqslant
e^{-t}\sum_{n=0}^{\infty}\frac{(\theta t)^n}{n!}\|b\|
=
e^{-t(1-\theta)}\|b\|.
\]
This proves \upshape(ii). Since the restriction of $(\T_t)$ to $\B$ is the identity semigroup, it has the $e$-property on $\B_{\mathrm{sa}}$. Theorem~\ref{thm:expected-subalgebra-lifting} now applies with $\Ldom_\B=\B_{\mathrm{sa}}$ and $\Ldom=\A_{\mathrm{sa}}$, proving \upshape(iii).
\end{proof}

\begin{corollary}\label{cor:compact-group-averaging}
Let $G$ be a compact group acting strongly continuously on a unital $C^*$-algebra $\A$ by $*$-automorphisms $\alpha_g$. Assume that the fixed-point algebra
\[
\B:=\A^G:=\{a\in\A:\alpha_g(a)=a\ \text{for all }g\in G\}
\]
admits the canonical conditional expectation
\[
E(a):=\int_G \alpha_g(a)\,dg.
\]
Let $\mu$ be a Borel probability measure on $G$, define the unital completely positive map
\[
P(a):=\int_G \alpha_g(a)\,d\mu(g),
\qquad a\in\A,
\]
and assume that there exists $\theta\in[0,1)$ such that
\[
\|P(b)\|\leqslant \theta\,\|b\|,
\qquad b\in\ker E.
\]
Then the semigroup generated by $P-\id$ satisfies the conclusion of Theorem~\ref{thm:contractive-expected-subalgebra}. If, in addition, $\B$ is a Cartan masa, this gives a Cartan-masa instance of the expected-subalgebra lifting mechanism.
\end{corollary}

\begin{proof}
The map $E$ is a conditional expectation onto $\B$, and Haar invariance gives $PE=EP=E$. Moreover, $P|_\B=\id_\B$ because every element of $\B$ is fixed by every $\alpha_g$. The conclusion therefore follows immediately from Theorem~\ref{thm:contractive-expected-subalgebra}.
\end{proof}

\begin{example}\label{ex:diagonal-lindblad}
Let $\A=M_n(\C)$, let $\D\subseteq M_n(\C)$ be the diagonal masa in the standard matrix units $(e_{ij})_{1\leqslant i,j\leqslant n}$, and let $E_\D$ be the diagonal conditional expectation. Choose self-adjoint diagonal matrices
\[
H=\operatorname{diag}(h_1,\dots,h_n),
\qquad
V_r=\operatorname{diag}(v_{r,1},\dots,v_{r,n})
\quad (1\leqslant r\leqslant m),
\]
and let $(\T_t)_{t\geqslant 0}$ be the uniformly continuous quantum Markov semigroup with generator
\[
\mathcal G(a)=i[H,a]+\sum_{r=1}^m\left(V_r aV_r-\frac12(V_r^2a+aV_r^2)\right).
\]
If
\[
\gamma:=\frac12\min_{i\neq j}\sum_{r=1}^m (v_{r,i}-v_{r,j})^2>0,
\]
then
\[
\|\T_t(a-E_\D(a))\|\leqslant \sqrt n\,e^{-\gamma t}\,\|a-E_\D(a)\|,
\qquad a\in M_n(\C),
\]
and consequently $(\T_t)$ has the $e$-property on $M_n(\C)_{\mathrm{sa}}$.
\end{example}

\begin{proof}
For the matrix units one computes
\[
\mathcal G(e_{ij})=
\left(i(h_i-h_j)-\frac12\sum_{r=1}^m(v_{r,i}-v_{r,j})^2\right)e_{ij}.
\]
Hence
\[
\T_t(e_{ij})=
\exp\!\left(t\,i(h_i-h_j)-\frac t2\sum_{r=1}^m(v_{r,i}-v_{r,j})^2\right)e_{ij}.
\]
In particular, $\T_t(e_{ii})=e_{ii}$, so $\T_t$ fixes $\D$ pointwise and commutes with $E_\D$. If
\[
b=a-E_\D(a)=\sum_{i\neq j}a_{ij}e_{ij},
\]
then
\[
\T_t(b)=\sum_{i\neq j} c_{ij}(t)a_{ij}e_{ij},
\qquad |c_{ij}(t)|\leqslant e^{-\gamma t}.
\]
Therefore the Hilbert--Schmidt norm satisfies
\[
\|\T_t(b)\|_2^2=
\sum_{i\neq j}|c_{ij}(t)|^2|a_{ij}|^2
\leqslant e^{-2\gamma t}\sum_{i\neq j}|a_{ij}|^2
=e^{-2\gamma t}\|b\|_2^2.
\]
Since $\|x\|\leqslant \|x\|_2\leqslant \sqrt n\,\|x\|$ on $M_n(\C)$, we get
\[
\|\T_t(b)\|\leqslant \|\T_t(b)\|_2\leqslant e^{-\gamma t}\|b\|_2\leqslant \sqrt n\,e^{-\gamma t}\|b\|.
\]
The conclusion now follows from Theorem~\ref{thm:expected-subalgebra-lifting}. This example is the simplest Cartan-masa instance of a dephasing Lindbladian, compare the qualitative finite-dimensional analyses of decoherence in \cite{BaumgartnerNarnhoferThirring2008,Rebolledo2005,CarboneSassoUmanita2014,CarboneSassoUmanita2015}.
\end{proof}

\section{Further questions}

\begin{problem}
Determine when approximation properties such as nuclearity, injectivity or finite nuclear dimension can be combined with positivity-improving estimates, approximate equivariance or functional inequalities to produce an exact Doeblin order lower bound.  Approximate finite-dimensional factorisations alone do not yield the required inequality \cite{BrownOzawa2008,CarlenMaas2017,KastoryanoTemme2013,KastoryanoBrandao2016}.
\end{problem}

\begin{problem}
While Theorem~\ref{thm:collapse} and Proposition~\ref{prop:jamison} settle the full-state-space hierarchy among the ordinary, eventual, Ces\`aro and pointwise eventual equicontinuity properties, it remains open to classify when the Ces\`aro $e$-property upgrades to the ordinary $e$-property under additional spectral or peripheral-point-spectrum assumptions on the generator.
\end{problem}

\begin{problem}
Find broad infinite-dimensional hypotheses under which decay in an $L^2$- or Dirichlet-form norm upgrades to the operator-norm decay required by the expected-subalgebra transfer theorem.
\end{problem}

\section*{Use of artificial intelligence}

This paper was prepared in 2024--2025 with some incremental additions in 2026. ChatGPT was used to polish the presentation. The authors take full responsibility for the content of the paper.

\end{document}